\documentclass[]{article}

\usepackage{graphicx} % for pdf, bitmapped graphics files
\usepackage[a4paper,left=35mm,right=35mm,top=30mm,bottom=30mm,marginpar=25mm]{geometry}
\usepackage{amsmath}
\usepackage{amsthm}
\usepackage{amssymb}
\usepackage{xcolor} %[monochrome] 
\usepackage[pdfpagelabels=false,colorlinks=true,linkcolor=blue,citecolor=red,urlcolor=cyan]{hyperref}
\usepackage{subcaption}
\usepackage[displaymath,mathlines]{lineno}
\allowdisplaybreaks
\newcommand{\R}{{\mathbb{R}}}

\newtheorem{definition}{Definition}
\usepackage{graphicx} % for pdf, bitmapped graphics files
\usepackage{amsmath}
\usepackage{amssymb}
\usepackage{xcolor} %[monochrome] 
\usepackage{hyperref}
\usepackage{tikz, pgfplots}
\usetikzlibrary{arrows,automata}

\allowdisplaybreaks

\newtheorem*{remark}{Remark}

\newcommand{\pddt}{\frac{\partial}{\partial t}}
\newcommand{\pddx}{\frac{\partial}{\partial x}}

\newcommand{\sub}[1]{_{#1}}%{^{(#1)}} % placement of the index of a
\newcommand{\subq}[1]{_{#1}} % placement of the index of q
\newcommand{\netg}{g} % Function g, boundary of Lambda

\newcommand{\fmax}{f^{\textnormal{max}}}

\newcommand{\rhozeronetwork}{\rho^0}
\newcommand{\rhoact}{\bar{\rho}}
\newcommand{\rhonum}{\tilde{\rho}}

\newcommand{\gridcell}{C}

\newcommand{\rhomaxnet}{\rho^{\textnormal{max}}}

\newcommand{\iin}{i} % Incoming arcs
\newcommand{\iout}{\hat{i}} % Outgoing arcs

\newcommand{\narc}{m_i} % Vertex y
\newcommand{\narcin}{m_{\iin}} % Vertex y
\newcommand{\smoothpar}{\xi} % Smoothing parameter

\newcommand{\hin}{h^{n,\text{in}}} % Numerical inflow
\newcommand{\hout}{h^{n,\text{out}}} % Numerical outflow
\newcommand{\intanf}{\bar{\mathfrak{a}}} % Numerical inflow
\newcommand{\intend}{\bar{\mathfrak{b}}} % Numerical outflow

\makeindex

\author{Adriano Festa \thanks{Institut National des sciences appliqu\'ees (INSA) Rouen, Laboratoire de Math\'ematiques, 685 Avenue de l'Universit\'e, 76800 Saint-\'Etienne-du-Rouvray. 
                {\tt\small adriano.festa@insa-rouen.fr} }
 \and Simone G\"ottlich
               \thanks{University of Mannheim, Department of Mathematics, A5-6, 68131 Mannheim, Germany.
                {\tt\small goettlich@uni-mannheim.de} }
                \and Marion Pfirsching\thanks{University of Mannheim, Department of Mathematics, A5-6, 68131 Mannheim, Germany. 
                {\tt\small mpfirsch@mail.uni-mannheim.de} }
}

\title{\Large \bf A model for a network of conveyor belts with discontinuous speed and capacity \thanks{This work was partially supported by the Haute-Normandie Regional Council via the M2NUM project and the project GO 1920/7-1
by the German Research Foundation (DFG).}}

\begin{document}

\maketitle
%\date{}

\begin{abstract}
We introduce a macroscopic model for a network of conveyor belts with various speeds and capacities. In a different way from traffic flow models, 
the product densities are forced to move with a constant velocity unless they reach a maximal capacity and start to queue. 
This kind of dynamics is governed by scalar conservation laws consisting of a discontinuous flux function.
We define appropriate coupling conditions to get well-posed solutions at intersections and provide a detailed description of the solution.
Some numerical simulations are presented to illustrate and confirm the theoretical results for different network configurations.
\end{abstract}

%% REQUIRED
%\begin{keywords}
%    production systems, conservation laws with discontinuous flux, numerical simulations
%\end{keywords}
%
%% REQUIRED
%\begin{AMS}
%  90B30, 35L65, 65M25, 
%\end{AMS}

\section{Introduction}

Conveyor belts have attained a favored position in transporting bulk materials due to their economy, reliability, safety, versatility and almost unlimited range of variations, conveying a wide variety of materials. 
However, theoretical study of these systems is non-trivial since the presence of intersections, fixed or smart divert/merge devices, and different production speeds add complexity to the structure and generate interesting dynamics.
For a general presentation of the models in the literature, we refer to the monograph \cite{d2010modeling} and the references therein.
As more recent contributions, we mention in particular the works \cite{armbruster2011scalar,Garavello2007,gottlich2013discontinuous,herty2013existence,Wiens2013}, where the authors introduce a production (or traffic) model with discontinuous flux and study the properties of solutions. 
While in \cite{armbruster2011scalar,gottlich2013discontinuous,Wiens2013}, the focus is on the mathematical modeling, and numerical simulation of the discontinuous flux function, the authors in \cite{Garavello2007,herty2013existence} prove the existence of solutions using wave-front tracking.  

In this work, we consider a production network consisting of multiple linked conveyor belts. 
We remark that the new network model differs essentially from \cite{gottlich2013discontinuous} in the proposed coupling conditions and the concept of solutions.
Each arc in the network corresponds to a conveyor belt with a certain constant speed and capacity along the arc. We suppose that the capacity is limited meaning that a maximum value of the density cannot be exceeded. Since we assume no buffers at the intersections, in the case of capacity drop the products get stuck on the belt (but they are still transported with a certain velocity). This is a key difference to traffic flow models \cite{garavello2016models,garavello2006traffic},
where the velocity is dependent on the density and vehicles have to stop once a maximal capacity is reached.

\begin{figure}[!t]
    \begin{center}
        \includegraphics[width=7cm]{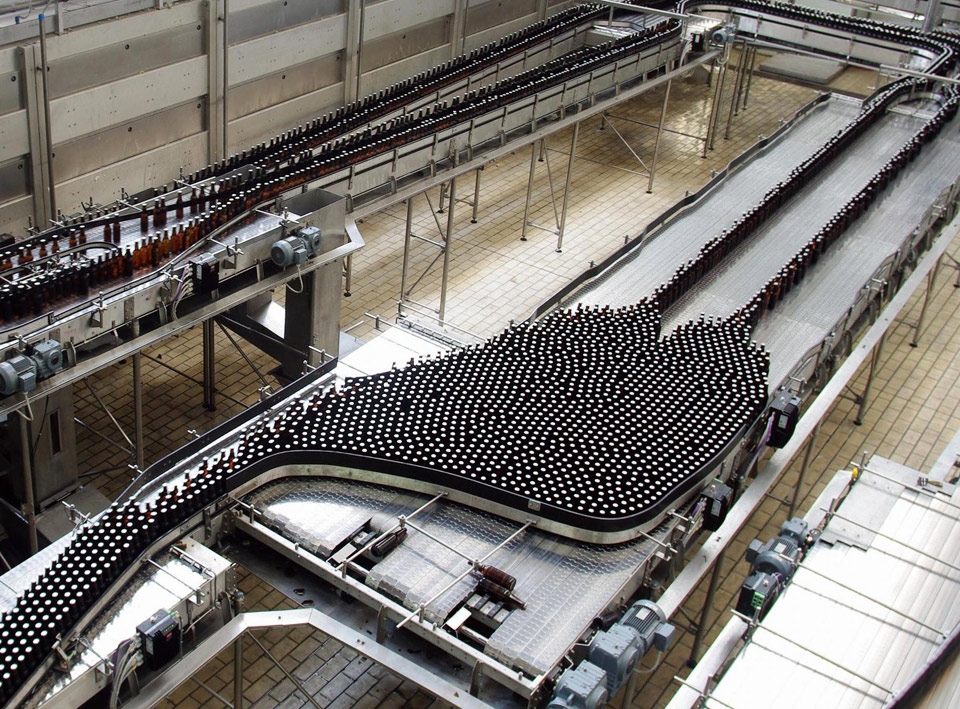}
        \caption{A conveyor belt in a brewery. \emph{Image courtesy of Sidel Blowing $\&$ Services SAS.}}
        \label{pic}
    \end{center}
\end{figure}

From an application viewpoint, the network model is particularly appropriate to study production systems for bottling, canning, and packaging. Figure \ref{pic} shows a brewery, where beer bottles are transported through a system of conveyor belts. 
We observe different lanes converging in an accumulator device and some ample space in the middle to stock the units that cannot be absorbed by the system.

The paper is organized as follows: in Section \ref{s:mod} we discuss the basic model and the notion of weak solutions which permits to derive the existence of a solution. We extend the framework to networks in Section \ref{s:net}: for one-to-one junction (Sec.~\ref{s121}) as well as diverging and merging intersection (Sect. \ref{s221} and \ref{s122}) where we build an analytic solution for the model. 
In Section \ref{snum}, we introduce a suitable discretization method to tackle the network model and present in Section \ref{stests} numerical results for different network settings where we show how the numerical approximations converge to the (previously described) analytic solution of interest.

\section{Basic model}\label{s:mod}

We start recalling the model originally introduced in \cite{armbruster2011scalar} and generalized to networks in \cite{gottlich2013discontinuous}. In both articles, a production system is described by a conservation law with discontinuous flux function and, if different from zero, constant speed $a>0$. More precisely for a closed set $\Omega\subset \R$ and calling $\rho:\Omega\times[0,T]\rightarrow [0,\rhomaxnet]$ the product density, the evolution of the system is
      \begin{equation}\label{eq:1}
\begin{cases}
\displaystyle  \partial_t\rho(x,t)+\partial_x \left(a H(\rhomaxnet-\rho(x,t))\rho(x,t)\right)=0,\\
\displaystyle \rho(x,0)=\rhozeronetwork(x),  
\end{cases}
\end{equation} 
where $H$ is the Heaviside function and the initial data $\rhozeronetwork$ is a function of bounded variation
satisfying $\rhozeronetwork(x)\leq \rhomaxnet$.

We point out that on a single arc (neglecting some possible boundary effects and assuming \eqref{eq:1} for all $x\in\Omega$) the equation simply reduces to an \emph{advection equation} and therefore the solution is simply 
\begin{equation}\label{sol:tran}
 \rho(x,t)=\rhozeronetwork(x-at).
 \end{equation}
This clearly differs from traffic models where the speed typically depends on the local density and shocks can appear on a single arc even with smooth initial data \cite{garavello2016models,garavello2006traffic}. The presence of discontinuities in the flux function, anyway, poses some problems in the case of $\rhozeronetwork(x)= \rhomaxnet$ for some $x\in\Omega$. 

A classic approach to deal with this problem is to use a standard regularization of the flux function using some Friedrichs' mollifiers. As it has been shown in \cite{dias2005solutions}, the solutions of the regularized problem converge to a \emph{bounded entropy weak solution}, in the particular sense of the definitions below. Therefore, we denote by $\tilde f$ the following multivalued function
$$\tilde f(\rho)=a \rho\quad  \hbox{ if }\rho\neq \rhomaxnet, \quad \tilde f(\rhomaxnet)=[0,a\rhomaxnet].$$

\begin{definition}\label{def:weak}
A function $\rho\in L^\infty(\R\times [0,T])$ is called \emph{weak solution} to the Cauchy problem \eqref{eq:1} if there exists a function $v\in L^{\infty}(\R\times [0,T])$ such that $v(x,t)\in\tilde f(\rho)$ a.e. and 
$$\int_0^T\int_\Omega \rho\frac{\partial\phi}{\partial t}dx dt+\int_0^T\int_\Omega v\frac{\partial\phi}{\partial x}dx dt+\int_\Omega \rhozeronetwork(x)\phi(x,0)dx=0$$
for each $\phi\in C^1_c(\R\times [0,T])$ (where $\phi\in C^1_c$ means $\phi\in C^1$ with compact support).
\end{definition}

Classically, the definition above is completed by the following notion of \emph{entropy weak solutions}:  denote by $\tilde H$ the following multivalued function
$$\tilde H(\rho)=H(\rho)\quad  \hbox{ if }\rho\neq 0, \quad \tilde H(0)=[0,1].$$

\begin{definition}\label{def:entropy}
A weak solution $\rho$ of the Cauchy problem \eqref{eq:1} is called an \emph{entropy weak solution} if, for each entropy $\eta\in C^1(\R)$, $\eta$ convex, there exists a function $w\in L^{\infty}(\R\times [0,T])$ such that $w(x,t)\in\tilde H(\rho(x,t))$ a.e. and 
$$\frac{\partial}{\partial t}\eta(\rho)+\frac{\partial}{\partial x}F(\rho)-\eta'(\rhomaxnet)\frac{\partial w}{\partial x}\leq 0, $$
where
$$F(\rho)=a\int_0^\rho \eta'(s)H(\rhomaxnet-s)ds.$$
\end{definition}

\begin{remark}
We observe that that the solution \eqref{sol:tran} is a weak entropy solution in the sense of Definitions \ref{def:weak} and \ref{def:entropy}. 
This can be shown by choosing
$$ v(x,t)=a\rho\in \tilde f(\rho), \quad \hbox{for }\rho\in [0,\rhomaxnet]$$
and
$$ \eta(\rho)=|\rho-k|, \quad \hbox{for }\rho\in [0,\rhomaxnet],\quad w(x,y)\equiv 1$$
for any constant $k\in\R$.

At the same time, as it is possible to see adapting the example provided in \cite{dias2004riemann}, in case of congested initial solution, a collection of weak solutions are acceptable as long as a condition (derived by the jump on the flux) on the speed of the congested area is verified. This means that in our case, for equation \label{eq:1}  with initial condition equal to 
$$ \rho^0(x)=\rho^{max}\,\chi_{(x<0)}, \quad x\in\R$$
there are at least two distinct entropy solutions equal to
$$ \widehat\rho(t,x)=\rho^{max}\,\chi_{(x/t<\alpha)}, \hbox{ and } \quad \overline\rho(t,x)=\chi_{(x/t<0)}\equiv\rho^0(x),\quad x\in\R.$$
%In other words, an area of maximal density may not move or move with a speed equal to the speed of the belt (advection formula). \\
Since we are interested in the more meaningful solution from the point of view of the  application, it seems clear that the advection formula provides the most meaningful candidate in the case of the presence of a congested region (differently from other applications e.g.  traffic models). This point has been discussed in detail in \cite{herty2013existence}, where the authors introduce the additional concept of congested/non congested region  in order to select such unique solution. Instead, in the present paper, we focus on building an entropy solution in the more general case of networks and showing how a suitable numerical scheme provides a good approximation of it.
\end{remark}

\section{Extension to networks} \label{s:net}

We now extend the model to a network represented by a directed graph $\Gamma=(V, E)$ with $V \neq \emptyset$ the set of vertices and $E\subset V\times V$ the set of arcs. For a fixed node $v \in V$, the sets $\delta_v^-, \delta_v^+$ denote the ingoing and outgoing arcs, respectively. 
We consider the following network problem:% on a network $\mathcal N:=\bigcup_{i\in E} \Omega_i$:
\begin{equation}
	\label{eq:model_network}
	\begin{cases}
	  \pddt \rho(x,t)+\pddx \left(f_i(\rho(x,t))\right)=0,  &x\in \Omega_i, ~ t \in [0, T]\\
	 \rho(x,0)=\rhozeronetwork_i(x),  &   x\in \Omega_i\\
	\sum\limits_{i \in \delta_v^-} f_i(\rho(v,t))= \sum\limits_{i\in \delta_v^+} f_i(\rho(v,t)),&  v\in V,
	\end{cases}
\end{equation} 
where
\begin{equation}
\label{eq:network_flux_arci}
f_i(\rho(x,t))=a\sub{i} H(\rhomaxnet_i-\rho(x,t))\rho(x,t), \quad x \in \overline{\Omega}_i, t \in [0,T]
\end{equation}
is the flux function on arc $i \in E$ and $a\sub{i} \in \R^+$ the transport velocity. 
At intersections, we assume the conservation of flux and 
to obtain well-posed solutions we need additional conditions which are discussed later in this section.
We underline that, differently from \cite{gottlich2013discontinuous}, speed and capacity are different for each arc
leading to new  discontinuities at the intersection points. 

\subsection{One-to-one junction}\label{s121}
The first case we consider is a one-to-one junction. It can be seen as a special case of a one-dimensional problem, as described by \eqref{eq:1}, with discontinuities in the velocity $a$ and the capacity $\rhomaxnet$.

For simplicity in the notation, we consider the problem on 
$\Omega =\Omega_1\cup \{0\}\cup \Omega_2=(-\infty,0)\cup \{0\}\cup(0,\infty),$
where the intersection is located at $ x= 0$. 
The model equations are given by \eqref{eq:model_network}, where the junction condition is specified as
\begin{equation}
%	\label{eq:network_1to1_model}
%	\begin{align}
% & \label{eq:network_1to1_transport} \pddt \rho(x,t)+\pddx \left(f_i(\rho(x,t))\right)=0, \quad (x,t)\in \Omega_i\times(0,T],\;i=1,2\\
% &\rho(x,0)=\rhozeronetwork_i(x), \qquad \qquad\quad\quad \quad ~ ~ ~ ~x\in \Omega_i,\;i=1,2\\
\label{eq:network_1to1_junctioncond} f_1(\rho(0,t))=f_2(\rho(0,t))
%\end{align} 
\end{equation} 
with flux function \eqref{eq:network_flux_arci} on arc $i = 1,2$.
The solution can be easily computed if no congestion occurs during the transportation, i.e.
\begin{equation}
\label{regular} a\sub{1}\rhozeronetwork_1(x)\leq a\sub{2}\rhomaxnet_2 \quad \forall x\in\Omega_1.
\end{equation}

In this case, we can derive the solution as follows.
For $x \in \Omega_1$, the characteristics of the problem are simply the straight lines $y(t)=x-a\sub{1} t$, which lead to the solution
\[\rho(x,t) = \rhozeronetwork_1(x-a\sub{1} t) \quad \text{for } (x,t): x<0.\]
%This corresponds to the solution of the advection equation.
Analogously, we find
\[\rho(x,t) = \rhozeronetwork_2(x-a\sub{2} t) \quad \text{for } \left\lbrace(x,t) \in \Omega_2\times (0,T] ~\big\lvert ~ x>a\sub{2} t \right\rbrace\]
which corresponds to the density initially placed on $\Omega_2$. The solution for $(x,t)$ in the case $0<x<a\sub{2} t$ is found by the junction condition \eqref{eq:network_1to1_junctioncond}:
$f_1(\rho(x^+,t))=f_2(\rho(x^-,t))$
where $x^\pm$ denotes $\lim_{h \rightarrow 0} x \pm h$. This leads to \[\rho(x^-,t)=\frac{a\sub{1}}{a\sub{2}}\rho(x^+,t)\]
at the interface $x=0$, which gives the solution of the intermediate part with adapted transport velocity.
The solution is then described by
\begin{equation}
\label{sol:reg}
\rho(x,t)=\begin{cases}
\rhozeronetwork_1(x-a\sub{1} t), & (x,t)\in \Omega_1\times(0,T]\\
\frac{a\sub{1}}{a\sub{2}}\rhozeronetwork_1\left(-a\sub{1}\left(t-\frac{x}{a\sub{2}}\right)\right), & (x,t) \in \Omega_2\times (0,T]: \;0 < x \leq a\sub{2} t\\
\rhozeronetwork_2(x-a\sub{2} t), & (x,t) \in \Omega_2\times (0,T]: \;x>a\sub{2} t.
\end{cases}
\end{equation}
%where $\overline{\Omega}_1$ denotes the closure of $\Omega_1$.

\begin{figure}[!t]
	\begin{center}
		\begin{tikzpicture}	
		\put(0,0){\node[draw=none]{\includegraphics[width=0.4\textwidth]{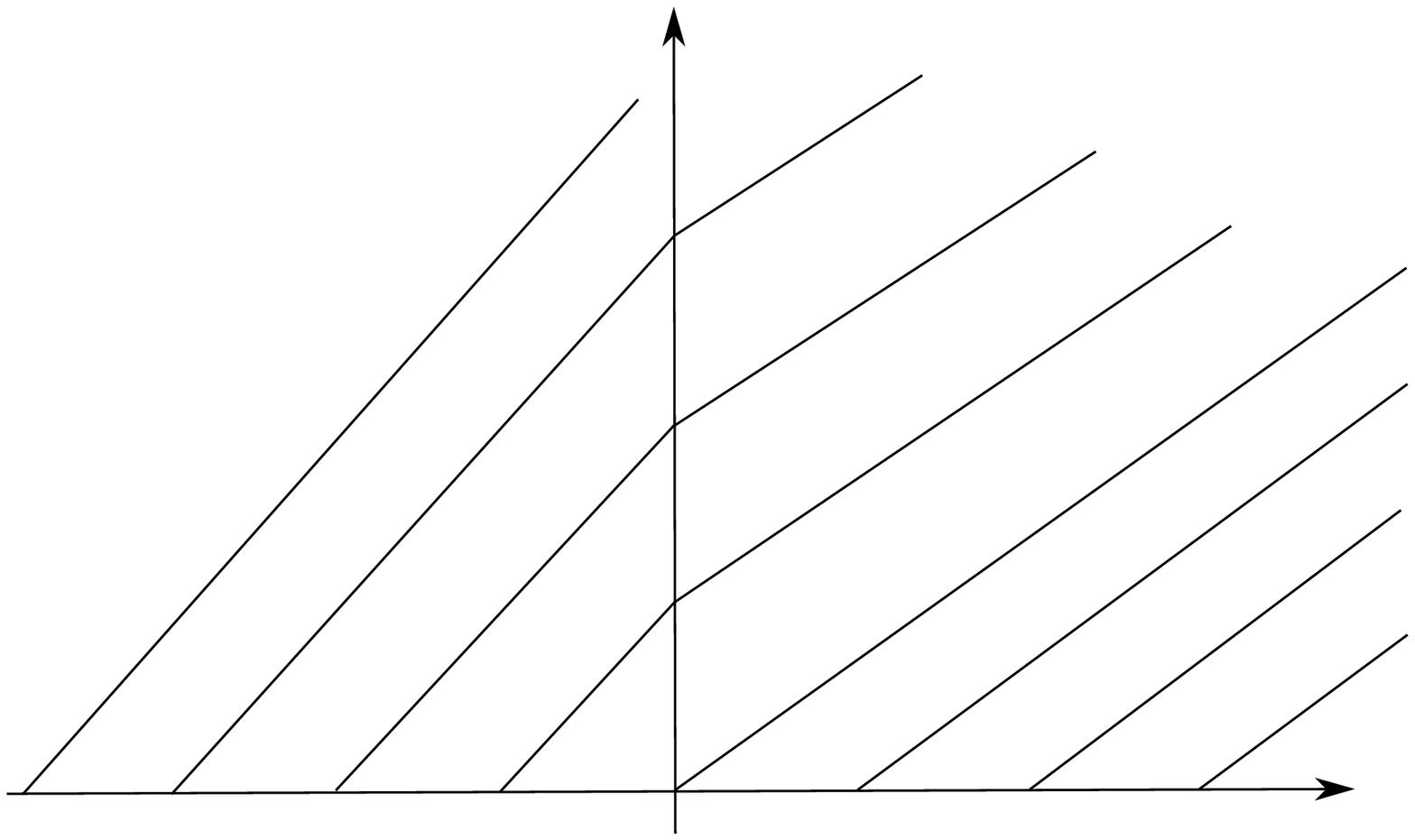}};}
		\put(0,0){\node[draw=none] at (-0.15, -1.7) {0};}
		\put(0,0){\node[draw=none] at (2.6, -1.5) {$x$};}
		\put(0,0){\node[draw=none] at (-0.15, 1.8) {$t$};}
		\put(0,0){\node[draw=none] at (-2.5, 0) {$x=a\sub{1}t+c$};}
		\put(0,0){\node[draw=none] at (3, 0.8) {$x=a\sub{2}t+c$};}
		\end{tikzpicture}
		\vspace{-0.3cm}
		\caption{Characteristics in the non-congested case}
		\label{probchar}
	\end{center}
\end{figure}

The characteristics of this solution are shown in Figure \ref{probchar} for $a\sub{1} < a\sub{2}$. 
%It can be interpreted as the trajectory a transported object would follow in the space-time diagram.
Note that the solution might be discontinuous at $x=0$ and $x = a\sub{2}t$. If the initial data $\rhozeronetwork_1$ and $\rhozeronetwork_2$ is continuous, 
the solution $\rho(x,t)$ keeps this property in all other points. Along each characteristic $c \in \R$ is constant.

\begin{remark}
	The interpretation of condition \eqref{regular} is as follows: We know that the flux is $a\sub{1} \rho(x,t)$ if the maximal density $\rhomaxnet_2$ is not reached. In this case, the capacity of arc 1 before the intersection has no influence. This is due to the positive velocities $a\sub{1}$ and $a\sub{2}$.
As we have already discovered, congested areas never appear in the interior of an arc, but only in conjunction with intersections.
\end{remark}
If condition \eqref{regular} is not satisfied, congestion arises and the problem becomes more involved. In particular, it is non-trivial to obtain a correct weak entropy solution using only Definitions \ref{def:weak} and \ref{def:entropy}. 
\begin{figure}[!t]
	\begin{center}	
		\begin{tikzpicture}	
		\put(0,0){\node[draw=none]{\includegraphics[width=0.5\textwidth]{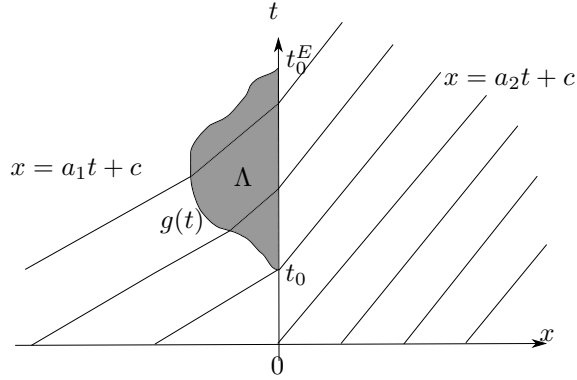}};}
		\put(0,0){\node[draw=none] at (-0.05, -2.3) {0};}
		\put(0,0){\node[draw=none] at (3.5, -1.9) {$x$};}
		\put(0,0){\node[draw=none] at (-0.1, 2.4) {$t$};}
		\put(0,0){\node[draw=none] at (-2.7, 0.3) {$x=a\sub{1}t+c$};}
		\put(0,0){\node[draw=none] at (3, 1.5) {$x=a\sub{2}t+c$};}
		\put(0,0){\node[draw=none] at (0.2, -1.1) {$t_0$};}
		\put(0,0){\node[draw=none] at (-1.3, -0.4) {$g(t)$};}
		\put(0,0){\node[draw=none] at (-0.5, 0.2) {$\Lambda$};}
		\put(0,0){\node[draw=none] at (0.25, 1.75) {$t_0^E$};}
		\end{tikzpicture}
		\vspace{-0.3cm}
		\caption{Trajectories in the congested case}
		\label{probchardis}
	\end{center}
\end{figure}

Let $t_0$ denote the first point, where condition \eqref{regular} is violated, i.e.~the first time of congestion:
\begin{equation}\label{t0}
t_0=\inf\left\{t\geq 0 \hbox{ such that }\rhozeronetwork_1(-a\sub{1} t)>\frac{a\sub{2}}{a\sub{1}}\rhomaxnet_2\right\}.
\end{equation} 
We track the interface describing the congested area at maximal density $\rhomaxnet_1$ that can appear in $\Omega_1$ and call the congested region $\Lambda$, see Figure \ref{probchardis}. 
The interface is a time-dependent function $g(t)$. The evolution of $g(t)$, starting at time $t_0$, can be derived by integrating the difference between the fluxes entering and exiting the region $\Lambda$ as well as the current density. The entering flux at time $t$ is given by $a\sub{1}\rhozeronetwork_1(-a\sub{1}t)$, the exiting flux by $a\sub{2}\rhomaxnet_2$ since \eqref{regular} is violated and the maximal density on the outgoing arc is reached. The resulting density 
in the congested region is $\rhozeronetwork_1(y)$ for $\netg(t)-a\sub{1}t \leq y \leq -a\sub{1}t$. Summarizing, this leads to 
\begin{equation}
%\label{eq:network_1to1_evolutionLambda}
\int\limits_{t_0}^t \left(a\sub{1}\rhozeronetwork_1(-a\sub{1}s)-a\sub{2}\rhomaxnet_2\right)ds= \int\limits_{\netg(t)}^{0} (\rhomaxnet_1-\rhozeronetwork_1(y-a\sub{1} t))dy. 
\end{equation}

%which is equivalent to
%\begin{equation*}
%-\int\limits_{-a\sub{1}t_0}^{-a\sub{1}t} \left(\rhozeronetwork_1(y)\right)dy-(t-t_0)a\sub{2}\rhomaxnet_2= -\netg(t)\cdot \rhomaxnet_1 - \int\limits_{\netg(t)-a\sub{1}t}^{-a\sub{1} t} \rhozeronetwork_1(y)dy.
%\end{equation*}

Rearranging the terms, we can describe the congested region $\Lambda$ as
\begin{equation}
\label{eq:network_1to1_defLambda}
\Lambda:=\left\{(x,t)\in \overline{\Omega}_1\times [t_0,t_0^E], \hbox{ such that } g(t)\le x \le 0\right\}
\end{equation}
with interface $g$ defined as $g:[0, \infty) \rightarrow \R_{\leq 0}$, where $t$ is mapped to the solution of 
\begin{equation}
\label{eq:network_1to1_defg}
-x+(t-t_0)\frac{\rhomaxnet_2}{\rhomaxnet_1}a\sub{2}-\frac{1}{\rhomaxnet_1}\int \limits_{x-a\sub{1} t}^{-a\sub{1} t_0}\rhozeronetwork_1(s)ds=0,
\end{equation}
if this is negative, and to zero otherwise.\\
The final time of congestion is
$ t^E=\min\left\{t\geq t_0 \hbox{ such that }g(t)=0\right\}.$

Figure \ref{probchardis} shows the trajectories of the problem in the case of congestion for $a\sub{1} > a\sub{2}$. Outside the region $\Lambda$, they correspond to the characteristics. If the initial data $\rhozeronetwork_1$ and $\rhozeronetwork_2$ are smooth, the solution may still become discontinuous in the interface $x=g(t)$, in $x=0$ and in $x = a\sub{2}t$. The congested region $\Lambda$ is highlighted in gray.

It is straightforward to see that the whole congested region can be constituted by multiple non-connected sets, if the congestion disappears and condition \eqref{regular} is violated again. In that case, we set $t_0^E := t^E$ and there exists a $t_k$ such that
\[t_k=\inf\left\{t\geq t^E_{k-1} \hbox{ such that }\rhozeronetwork_1(-a\sub{1} t)>\frac{a\sub{2}}{a\sub{1}}\rhomaxnet_2\right\}, \quad k = 1, 2, \ldots ~.\]
The procedure to build this second connected subset of the whole congested region is the same as for the first one and so on. Therefore, it is sufficient to only consider one connected set $\Lambda$.

Inside the region $\Lambda$, the transport velocity $\bar{a}$ is such that the coupling condition \eqref{eq:network_1to1_junctioncond} holds true, i.e.~the inflow $\bar{a}\rhomaxnet_1$ equals the outflow $a\sub{2} \rhomaxnet_2$ at $x=0$. Therefore, the velocity inside the region $\Lambda$ is
\begin{equation}
\label{eq:network_1to1_bara}
\bar a=a\sub{2}\frac{\rhomaxnet_2}{\rhomaxnet_1}.
\end{equation}

\begin{remark}
	Since we assumed that condition \eqref{regular} is violated, we find a value $\bar x\in\Omega_1$ such that
	$a\sub{1}\rhozeronetwork_1(\bar x)>a\sub{2}\rhomaxnet_2 $ and $ \rhozeronetwork_1(\bar{x})\leq \rhomaxnet_1,$
	where the second term is due to the choice $\rho^0(x)\leq\rho^{\max}$. This implies
	\[a\sub{1}\rhomaxnet_1\geq a\sub{1} \rhozeronetwork_1(\bar x)>a\sub{2} \rhomaxnet_2.\]
	Dividing the first and the last term of this inequality by $\rhomaxnet_1$, we obtain 
	$\bar a<a\sub{1}$,
	which means that the velocity always decreases as soon as the mass enters the maximal density area $\Lambda$. This confirms the intuitive assumption that the transport velocity is reduced in the congested region $\Lambda$.
\end{remark}
Next, we derive a general solution for this case (shown in Figure \ref{probchardis}) for the modifications in the characteristics.

We follow an approach considering the associated Riemann problem. %, meaning that $\rhozeronetwork_1$ and $\rhozeronetwork_2$ are assumed to be constant.
We calculate the solution for a general initial condition, approximated by piecewise constant functions. The same approach is discussed in \cite{herty2013existence} for more general problems.

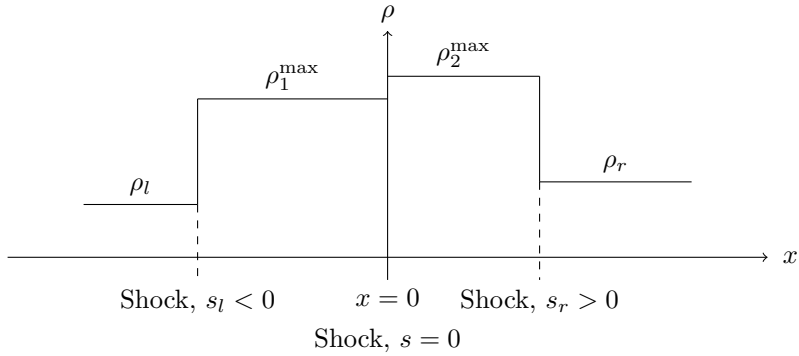
\begin{figure}[!t]
	\centering
	\begin{tikzpicture} 
	\draw[->] (-5,0) --  (5,0);
	\draw[->] (0,-0.3)node[below] {$x=0$} --  (0, 3) node[above] {$\rho$};
	\node[draw=none] at (0, -1.1)() {Shock, $s = 0$} ;
	\draw[-] (0,2.4) --  node[above] {$\rhomaxnet_2$} (2, 2.4) ;
	\draw[-] (2, 2.4) --  (2, 1) ;
	\draw[-] (2, 1)  --  node[above] {$\rho_r$} (4,1) ;
	\draw[dashed] (2, 1) -- (2, -0.3) node[below] {Shock, $s_r > 0$};
	\draw[-] (-2.5, 2.1) --  node[above] {$\rhomaxnet_1$} (0, 2.1) ;
	\draw[-] (-2.5, 2.1) --  (-2.5, 0.7) ;
	\draw[-] (-2.5, 0.7)  --  node[above] {$\rho_l$} (-4, 0.7) ;
	\draw[dashed] (-2.5, 0.7)  -- (-2.5, -0.3) node[below] {Shock, $s_l < 0$};
	\node[draw=none] at (5.3, 0)() {$x$} ;
	\end{tikzpicture}
	\caption{Solution in the congested case: evolution of three shock waves}
	\label{fig:shock_evolution}
\end{figure}

We briefly discuss the Riemann problem on $[0, T)\times \R$, consisting of equation \eqref{eq:model_network} and the initial data 
\begin{equation}
\rho(x,0)=\left\{\begin{array}{ll}
\rho_l, & x<0,\\
\rho_r, & x\geq 0.
\end{array}\right.
\end{equation}
If we allow waves with negative velocity, we can solve this Riemann problem by distinguishing the following cases:

\begin{itemize}
	\item[a)] $f_1 \leq \fmax_2$: 
	%$ a\sub{1}\rhozeronetwork_1(x)\leq a\sub{2}\rhomaxnet_2 $: 
	no congestion arises, \eqref{regular} is verified and the solution is \eqref{sol:reg} with piecewise constant initial data
	$\rhozeronetwork_1(x) = \rho_l$ and $\rhozeronetwork_2(x) = \rho_r$.
	\item[b)] $f_1 > \fmax_2$: 
	%$ a\sub{1}\rhozeronetwork_1(x)\ge a\sub{2}\rhomaxnet_2 $: 
	a congestion arises since \eqref{regular} is not verified. In this case, the solution consists of three shock waves (see Figure \ref{fig:shock_evolution}) starting at $x = 0$:
	\begin{itemize}
		\item[i)] At $x = 0$, a shock with velocity $s= 0$ arises, where the solution $\rho$ jumps from $\rhomaxnet_1$ to $\rhomaxnet_2$. For the special case $\rhomaxnet_1 = \rhomaxnet_2$, there is no jump in the solution at $x = 0$.
		\item[ii)] We obtain a left-going shock wave, where the density jumps from $\rho_l$ to $\rhomaxnet_1$ with negative velocity $s_l < 0$ (computed according to the Rankine-Hugoniot condition):
		\[ s_l = \frac{\bar{a}\rhomaxnet_1 - f_1(\rho_l)}{\rhomaxnet_1 - \rho_l} = \frac{\fmax_2 - f_1(\rho_l)}{\rhomaxnet_1 - \rho_l}  < 0.\]
		This shock wave describes the left boundary of the congested area $\Lambda$, meaning that it follows the function $g(t)$. Therefore, the transport right of this shock is with velocity $\bar{a}$ as defined in \eqref{eq:network_1to1_bara}. In the case $\rho_l \rightarrow \rhomaxnet_1 $, the shock velocity tends to $- \infty$.
		\item[iii)] We obtain a right-going shock wave, where the density jumps from $\rhomaxnet_2$ to $\rho_r$ with positive velocity $s_r = a\sub{2} > 0$ (computed according to the Rankine-Hugoniot condition):
		\[ s_r = \frac{f_2(\rho_r) - f_2(\rhomaxnet_2)}{\rho_r - \rhomaxnet_2} = \frac{a\sub{2}\rho_r - a\sub{2}\rhomaxnet_2}{\rho_r - \rhomaxnet_2} = a\sub{2} > 0.\]
	\end{itemize}
\end{itemize}
%	\end{itemize}
In the case $t_0 > 0$, we are in the non-congested case for $t < t_0$. At $t = t_0$, congestion starts and the shock waves appear.

In the general case with non-constant initial data, we obtain the following solution:
\begin{equation}\label{sol121}
\rho(x,t)=
\begin{cases}
\rhozeronetwork_1(x-a\sub{1} t), & (x,t)\in \Omega_1\setminus \Lambda \times (0, T]\\
\rhomaxnet_1, & (x,t)\in \Lambda \times (0, T]\\
\rhomaxnet_2, & (x,t) \in \Omega_2 \times (0,T]: \;x\leq a\sub{2} t,g(t-\frac{x}{a\sub{2}})\neq 0\\
\frac{a\sub{1}}{a\sub{2}}\rhozeronetwork_1\left(-a\sub{1}\left(t-\frac{x}{a\sub{2}}\right)\right), & (x,t)\in \Omega_2 \times (0,T]: \;x\leq a\sub{2} t, g(t-\frac{x}{a\sub{2}})= 0\\
\rhozeronetwork_2(x-a\sub{2} t), & (x,t)\in \Omega_2 \times (0,T]: \;x>a\sub{2} t.
\end{cases}
\end{equation} 
If $g(t)=0$ for all $t \in [0,T]$, we have $\Lambda = \emptyset$ and recover the previously described solution \eqref{sol:reg}.

Contrary to traffic flow models, rarefaction waves do not appear in the conveyor belt problem. This is due to the linear flux function up to the discontinuity and drops to zero.

\subsection{One-to-two junction}\label{s122}

\begin{figure}[htb]
    \centering
    \begin{minipage}{4cm}
%    \subfloat[Passive junction\label{Fig:122_passive}] {
        \begin{tikzpicture}    
        \put(0,0){\node[draw=none]{\includegraphics[width=0.8\textwidth]{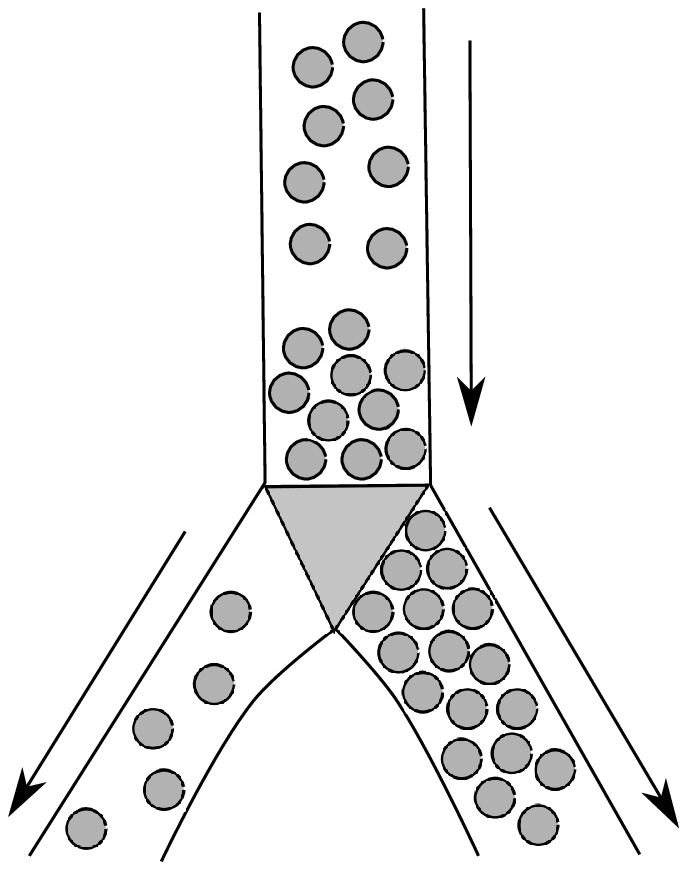}};}
        \put(0,0){\node[draw=none] at (1.0, 1.3) {$f_1$};}
        \put(0,0){\node[draw=none] at (-1.5, -1.0) {$f_2$};}
        \put(0,0){\node[draw=none] at (1.3, -0.8) {$f_3$};}
        \put(0,0){\node[draw=none] at (0,-0.5) {$D$};}
        \end{tikzpicture}%}
%        \caption{Passive junction}
%        \label{Fig:122_passive}
    \end{minipage}
    \begin{minipage}{4cm}
%    \subfloat[Active junction\label{Fig:122_active}] {
        \begin{tikzpicture}    
        \put(0,0){\node[draw=none]{\includegraphics[width=0.8\textwidth]{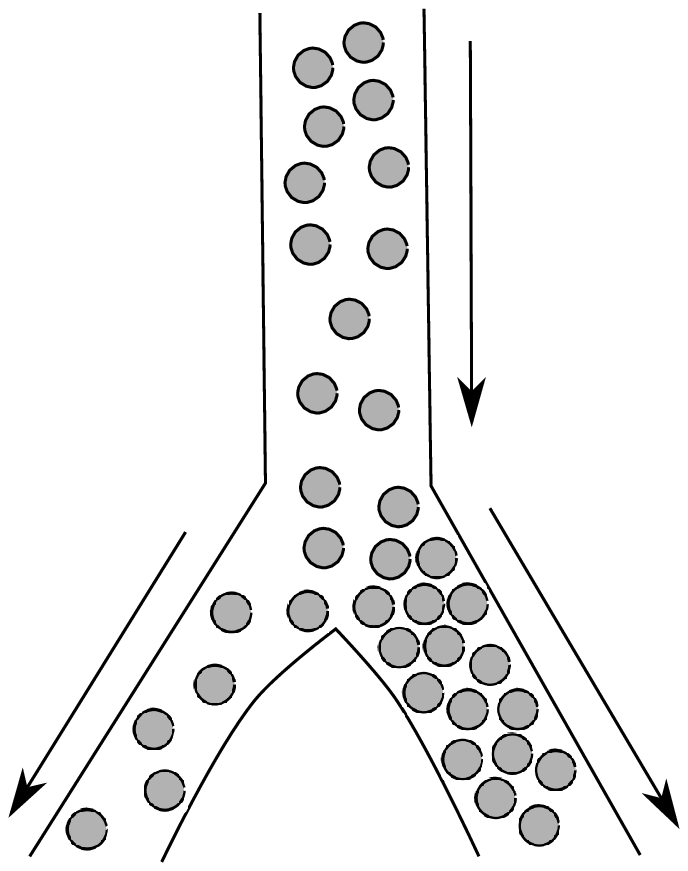}};}
        \put(0,0){\node[draw=none] at (1.0, 1.3) {$f_1$};}
        \put(0,0){\node[draw=none] at (-1.5, -1.0) {$f_2$};}
        \put(0,0){\node[draw=none] at (1.3, -0.8) {$f_3$};}
        \end{tikzpicture}%}
%        \caption{Active junction}
%        \label{Fig:122_active}
    \end{minipage}
    \caption{Scheme of the two cases considered of one-to-two junction: passive (left) and active (right)}
    \label{Fig:122}
\end{figure}

Now, we consider the case of a splitting intersection, where one arc is separated into two. We denote by $i = 1$ the incoming and by $i = 2,3$ the outgoing 
arcs (see Figure \ref{Fig:122}).

The choice of a distribution rule depends on the application:  if the intersection is ``passive'', we impose a fixed rate between the outgoing fluxes which is kept constant during the transportation process via a device $D$ (see the left picture in Figure \ref{Fig:122}). In this situation, congestion arises even if the outgoing arcs are not both congested. A different approach is an ``active'' junction (see right picture in Figure \ref{Fig:122}). In this case, the ratio between the outgoing fluxes can change at the intersection. The aim is the maximization of the total outgoing flux and the reduction of congestion. 
We underline that this behavior partially differs from the standard coupling conditions introduced for vehicular traffic fluxes (cf.~\cite{garavello2016models}), since the choice of the vehicles (left/right at a junction) cannot be determined by the local status of the traffic.

We consider the problem on 
$\Omega=\Omega_1\, e_1 \cup \{0\} \cup \Omega_2\, e_1\cup\Omega_3\, e_2=(-\infty,0)\,e_1\cup\{0\} \cup(0,\infty)\,e_1\cup(0,\infty)\,e_2,$
where $(e_1,e_2)$ is the standard base of $\R^2$, and we identify the element $x \in \Omega_1$ with the vector $(x,0)^T$, and analogously for arcs $i = 2,3$.
% state $x=x\,e_1$ if $x\in\Omega_1\cup \Omega_2$ and $x=x\,e_2$ if $x\in\Omega_3$. 
The equations that we consider are \eqref{eq:model_network}, where the junction condition is specified as
\begin{equation}
%    \label{eq:network_1to2_model}
%\begin{alignat}{2}
%  & \label{eq:network_1to2_transport} \pddt \rho(x,t)+\pddx \left(f_i(\rho(x,t))\right)=0, \qquad  \qquad&&(x,t)\in \Omega_i\times(0,T],\;i=1,2,3\\
% & \label{eq:network_1to2_initcond} \rho(x,0)=\rhozeronetwork_i(x), \qquad \qquad \qquad \qquad \qquad ~ ~ ~ ~ &&x\in \Omega_i,\;i=1,2,3\\
\label{eq:network_1to2_junctioncond} f_1(\rho(0,t))=  f_2(\rho(0,t))+ f_3(\rho(0,t)), 
%\end{alignat}
\end{equation} 
equipped with flux function \eqref{eq:network_flux_arci} on arc $i = 1, 2, 3$.

\medskip

\noindent{\bf Case 1: ``passive'' junction.} At first, we consider the case of same flux rates between the two exiting arcs. An intersection device $D$ keeps the ratio of the two outgoing fluxes constant, i.e., for a fixed distribution parameter $\mu\in[0,1]$, it holds
\begin{equation}
\label{eq:network_1to2_splittingcond}
f_2 = \mu f_1, \quad f_3 = (1-\mu)f_1,
\end{equation}
even if only one outgoing arc is congested.
%We highlight that condition \eqref{eq:network_1to2_splittingcond} verifies the conservation of fluxes at the junction point $x=0$.
The case $\mu \in \{0, 1\}$ reduces the problem to the one-to-one junction, so we consider $\mu \in (0,1)$ now.
The relation \eqref{eq:network_1to2_splittingcond} states that a fixed rate 
\begin{equation}
\label{eq:network_1to2_ratio}
f_2=\frac{\mu}{1-\mu}f_3
\end{equation}
is kept between the two outgoing fluxes during the evolution of the system, independent of the incoming flux. If no congestion arises, i.e.,
\begin{equation}\label{regular122} 
\rhozeronetwork_1(x) \leq \min\left\{\frac{1}{\mu}\frac{a\sub{2}}{a\sub{1}}\rhomaxnet_2,\frac{1}{1-\mu}\frac{a\sub{3}}{a\sub{1}}\rhomaxnet_3\right\}, \quad x\in\Omega_1
\end{equation}
holds true, the solution is obtained in a similar way as in Section \ref{s121}. We skip this point and draw our attention directly to a general formula for the solution (with or without congested areas).

Due to the constant rate \eqref{eq:network_1to2_ratio} between the two outgoing fluxes $f_2$ and $f_3$, it might happen that only one outgoing arc reaches the maximal density before congestion on the incoming arc $1$ arises. Therefore, we define the actual density on arc $2$ and $3$ as
\begin{equation*}
\rhoact_2=\min\left\{\frac{\mu}{1-\mu}\frac{a\sub{3}}{a\sub{2}}\rhomaxnet_3,\rhomaxnet_2\right\}, \quad \rhoact_3=\min\left\{\frac{1-\mu}{\mu}\frac{a\sub{2}}{a\sub{3}}\rhomaxnet_2,\rhomaxnet_3\right\}.
\end{equation*}
This is the minimum of the density the arc is supposed to take due to the distribution parameter, and the maximal density possible on this arc.
The first time of congestion $t_0$ can be determined as
\begin{equation*}
t_0=\inf\left\{t\geq 0 \hbox{ such that }\rhozeronetwork_1(-a\sub{1} t)>\min\left\{\frac{1}{\mu}\frac{a\sub{2}}{a\sub{1}}\rhomaxnet_2,\frac{1}{1-\mu}\frac{a\sub{3}}{a\sub{1}}\rhomaxnet_3\right\}\right\}.
\end{equation*}
We can introduce the interface $g(t)$, in analogy to \eqref{eq:network_1to1_defg}, with exiting flux $\rhoact_2 a\sub{2} + \rhoact_3 a\sub{3}$. The interface is defined as $g:[0, \infty) \rightarrow \R_{\leq 0}$, where $t$ is mapped to the solution of equation
\begin{equation}
\label{eq:network_1to2_passive_defg}
-x+(t-t_0)\frac{\rhomaxnet_2}{\rhomaxnet_1}a\sub{2}-\frac{1}{\rhomaxnet_1}\int\limits_{x-a\sub{1} t}^{-a\sub{1} t_0}\rhozeronetwork_1(s)ds=0,
\end{equation}
if this is negative, and to zero otherwise.
The region of congestion $\Lambda$ on $\overline{\Omega}_1$ is given by
\begin{equation}
\label{laambda}
\Lambda:=\left\{(x,t)\in \overline{\Omega}_1\times [t_0,t_0^E], \hbox{ such that } g(t)\leq x \leq 0\right\},
\end{equation}
analogously to \eqref{eq:network_1to1_defLambda}.
The general solution $\rho$ on $\Omega$ is given by
\begin{equation}\label{sol122_ratio}
\rho(x,t)=
\begin{cases}
\rhozeronetwork_1(x-a\sub{1} t), & (x,t)\in \Omega_1\setminus \Lambda \times(0,T]\\
\rhomaxnet_1, & (x,t)\in  \Lambda \times(0,T]\\
\rhoact_i, & (x,t)\in  \Omega_i \times(0,T]: \begin{aligned}[t] & x\leq a\sub{i} t, \\ & g(t-\frac{x}{a\sub{i}}) \neq 0 , i=2,3\end{aligned}\\
\alpha_i\frac{a\sub{1}}{a\sub{i}} \rho^{1}_{0}  \left(-a\sub{1}\left(t-\frac{x}{a\sub{i}}\right)\right), & (x,t) \in \Omega_i \times(0,T]: \begin{aligned}[t] & x\leq a\sub{i} t, \\ & g(t-\frac{x}{a\sub{i}}) = 0 , i=2,3\end{aligned}\\
%\;x \leq a\sub{i} t, g(t-\frac{x}{a\sub{i}})= 0, i=2,3\\
\rhozeronetwork_i(x-a\sub{i} t), & (x,t) \in \Omega_i \times(0,T]: \;x>a\sub{i} t, i = 2,3,
\end{cases}
\end{equation}
where
%\begin{equation}\label{alpha}
$\alpha_2=\mu \hbox{ and } \alpha_3=1-\mu$
%\end{equation}
for a compact notation.
\begin{remark}
    Congestion occurs if the maximal density of one of the two exiting arcs $i = 2,3$ is reached. 
    The other arc, even if the maximal density is not reached, shows a similar congested behavior, i.e., a value less than $\rhomaxnet_i$ is reached and kept. Since congestion arises without using the full capacity of both outgoing arcs, the duration of the congested phase is prolonged. 
\end{remark}

\noindent{\bf Case 2: ``active'' junction.}
We consider the possibility of a diverter interpreted as a device that keeps a constant ratio among the two outgoing fluxes as long as no arc is congested. If congestion arises, the device adapts the fluxes to ensure the maximal total outgoing flux, i.e., $f_2 = \fmax_2$ and $f_3 = \fmax_3$.

As before, we fix a parameter $\mu \in [0,1]$ in order to set a constant ratio in the non-congested case. As for the passive junction case \eqref{eq:network_1to2_splittingcond}, we have
$\mu f_1=f_2 \hbox{ and } (1-\mu)f_1=f_3.$
In order to get a unique solution also in the congested case, we define parameters $\beta_i$ corresponding to $\alpha_i$ in \eqref{sol122_ratio} with the following properties:
\begin{itemize}
	\item[a)] The flux conservation $f_1=f_2+f_3$ at the coupling is satisfied.
	\item[b)] The parameters $\beta_i$ are equal to $\alpha_i$, $i=2,3$ as in the previous case if no congestion occurs, i.e., condition \eqref{regular122} holds true.
	\item[c)] If only one outgoing arc $i$ is congested, i.e.,
	$\rhozeronetwork_1(x-a\sub{1} t)> \frac{a\sub{i}}{\alpha_i a\sub{1}}\rhomaxnet_i$
	at time $t$ the parameter $\beta_i$ changes to
	\[\beta_i = \frac{a\sub{i}}{a\sub{1}}\frac{\rhomaxnet_i}{\rhozeronetwork_1(-a\sub{1}t)}\]
	which is the smallest value to avoid congestion.% Hence, $\beta_i, ~ i = 2,3$ has to be adapted accordingly.
	\item[d)] A further change is necessary if the value 
	\[\beta_i = \frac{a\sub{i}\rhomaxnet_i}{a\sub{2}\rhomaxnet_2+a\sub{3}\rhomaxnet_3}\]
	is reached. This is the optimal ratio to maximize the flux through the junction. At this point the congestion starts.
\end{itemize}
Having stated the conditions on $\beta_i, ~ i = 2,3$, we can now derive the solution.
No congestion arises as long as the following inequality holds true:
\begin{equation}\label{regular122b} 
\rhozeronetwork_1(x) \leq \frac{a\sub{2}}{a\sub{1}}\rhomaxnet_2+\frac{a\sub{3}}{a\sub{1}}\rhomaxnet_3, \quad x\in\Omega_1.
\end{equation}
If condition \eqref{regular122b} is not fulfilled, the time $t_0$, meaning the first time of congestion, is independent of the distribution parameter $\mu$. It is then given by
\begin{equation*}
t_0=\inf\left\{t\geq 0 \hbox{ such that }\rhozeronetwork_1(-a\sub{1} t)>\frac{a\sub{2}}{a\sub{1}}\rhomaxnet_2+\frac{a\sub{3}}{a\sub{1}}\rhomaxnet_3\right\}.
\end{equation*}
In this case, the interface $g$ is defined as $g:[0, \infty) \rightarrow \R_{\leq 0}$, where $t$ is mapped to the solution of equation
\begin{equation}
-x+(t-t_0)\left(\frac{ \rhomaxnet_2}{\rhomaxnet_1}a\sub{2}+\frac{ \rhomaxnet_3}{\rhomaxnet_1}a\sub{3}\right)-\frac{1}{\rhomaxnet_1}\int\limits_{x-a\sub{1} t}^{-a\sub{1} t_0}\rhozeronetwork_1(s)ds=0,
\end{equation}
if this solution is negative, and to zero otherwise.
The region of congestion $\Lambda$ on $\Omega_1$ is given by \eqref{laambda}.
Contrary to the passive junction case, the function $g(t)$ considers the maximal flux on arcs $2$ and $3$ but no longer the distribution parameter $\mu$.

We define for $i = 2,3$
\begin{equation}\label{Beta}
\beta_i(t):=\min\left\{\max\left\{\alpha_i,\frac{a\sub{i}}{a\sub{1}}\frac{\rhomaxnet_i}{\rhozeronetwork_1(x-a\sub{1} t)}\right\}, \frac{a\sub{i}\rhomaxnet_i}{a\sub{2}\rhomaxnet_2+a\sub{3}\rhomaxnet_3}\right\}
\end{equation}
with $\alpha_2=\mu$ and $\alpha_3 = 1-\mu$. 
The general solution on $\Omega$ is then
\begin{equation}\label{sol122_active}
\rho(x,t)=
\begin{cases}
\rhozeronetwork_1(x-a\sub{1} t), & (x,t)\in \Omega\setminus \Lambda \times(0,T]\\
\rhomaxnet_1, & (x,t)\in  \Lambda \times(0,T]\\
\rhomaxnet_i, & (x,t) \in \Omega_i \times (0,T]: \begin{aligned}[t] & x\leq a\sub{i} t, \\ & g(t-\frac{x}{a\sub{i}}) \neq 0 \end{aligned}\\
\beta_i(t-\frac{x}{a\sub{1}})\frac{a\sub{1}}{a\sub{i}} \rho^{1}_{0}  \left(-a\sub{1}\left(t-\frac{x}{a\sub{i}}\right)\right), & (x,t)\in \Omega_i \times (0,T]:\begin{aligned}[t] & x\leq a\sub{i} t, \\ & g(t-\frac{x}{a\sub{i}}) = 0 \end{aligned}\\
% \;x \leq a\sub{i} t, g(t-\frac{x}{a\sub{i}}) = 0\\
\rhozeronetwork_i(x-a\sub{i} t), & (x,t)\in \Omega_i \times (0,T]: \;x>a\sub{i} t
\end{cases}
\end{equation}
with $i = 2,3$.
%We state the differences with the solution .
\begin{remark}
Compared to the ``passive'' junction \eqref{sol122_ratio}, congestion only occurs if the maximal capacity of both exiting arcs is reached. This implies that the choice of $\mu \in \{0, 1\}$ does no longer reduce to a one-to-one junction. 
%As a consequence, if we want to model an accumulation buffer, we can consider one exiting arc with a very high capacity, which is unused if no congestion arises. In the congested case, it will absorb the part of the flux that exceeds the capacity of the other arc.
\end{remark}

\subsection{Two-to-one junction}\label{s221}

In this part, we focus on the case of a merging junction. 
We know from traffic flow that in the free flow regime no additional information is needed. Conversely, in the congested case, we need a priority rule between the two incoming arcs, i.e., how to use released capacities of the outgoing arc. We denote by $i = 1,2$ the incoming and by $i = 3$ the outgoing arcs.

We consider the problem on 
$\Omega=\Omega_1\,e_1\cup \Omega_2\,e_2\cup\{0\}\cup\Omega_3\,e_1=(-\infty,0)\,e_1\cup (-\infty,0)\,e_2\cup\{0\}\cup(0,\infty)\,e_1$
with the same interpretation as before, where the system is given by \eqref{eq:model_network}, the coupling condition reads as
\begin{equation}
%\label{eq:network_2to1_model}
%\begin{alignat}{2}
%& \label{eq:network_2to1_transport} \pddt \rho(x,t)+\pddx \left(f_i(\rho(x,t))\right)=0, \qquad \qquad &&(x,t)\in \Omega_i\times(0,T],\;i=1,2,3\\
%& \label{eq:network_2to1_initcond} \rho(x,0)=\rhozeronetwork_i(x), &&x\in \Omega_i,\;i=1,2,3\\
\label{eq:network_2to1_junctioncond} f_1(\rho(0,t))+ f_2(\rho(0,t))=f_3(\rho(0,t)),
%\end{alignat}
\end{equation} 
and the flux function is again \eqref{eq:network_flux_arci} on arcs $i = 1, 2, 3$. % Equation \eqref{eq:network_2to1_transport} describes the transport on each arc, equation \eqref{eq:network_2to1_initcond} states the initial condition and condition \eqref{eq:network_2to1_junctioncond} describes the condition at the junction point, analogous to \eqref{eq:model_network}.\\
The solution can be directly computed if it holds
%\begin{equation}\label{regular221} 
%a\sub{1}\rhozeronetwork_1(x)+ a\sub{2}\rhozeronetwork_2(y)\leq a\sub{3}\rhomaxnet_3, \quad x\in\Omega_1, y\in\Omega_2, \hbox{ and }x=y,
%\end{equation}
\begin{equation}\label{regular221} 
a\sub{1}\rhozeronetwork_1(-a\sub{1}t)+ a\sub{2}\rhozeronetwork_2(-a\sub{2}t)\leq a\sub{3}\rhomaxnet_3,
\end{equation}
i.e.~no congestion arises. Then, the solution is given by
\begin{equation}\label{sol:reg221}
\rho(x,t)=
\begin{cases}
\rhozeronetwork_i(x-a\sub{i} t), & (x,t)\in \Omega_i\times(0,T],\, i = 1,2\\
%\rhozeronetwork_2(x-a\sub{2} t), & (x,t)\in \Omega_2\times(0,T],\\
\sum\limits_{i=1}^2 \frac{a\sub{i}}{a\sub{3}}\rhozeronetwork_i\left(-a\sub{i}\left(t-\frac{x}{a\sub{3}}\right)\right),
& (x,t)\in \Omega_3 \times (0,T]: x\leq a\sub{3} t \\
%\frac{a\sub{1}}{a\sub{3}}\rhozeronetwork_1\left(-a\sub{1}\left(t-\frac{x}{a\sub{3}}\right)\right)+\frac{a\sub{2}}{a\sub{3}}\rhozeronetwork_2\left(-a\sub{2}\left(t-\frac{x}{a\sub{3}}\right)\right),& \hspace{-0.2cm}(x,t)\in \overline{\Omega}_3 \times (0,T]: x\leq a\sub{3} t\\
\rhozeronetwork_3(x-a\sub{3} t), & (x,t) \in \Omega_3 \times (0,T]: x>a\sub{3} t
\end{cases}
\end{equation} 
and the incoming mass can be totally absorbed by the outgoing arc. This is independent of the capacity of the incoming arcs, and no further priority rule is needed to obtain a unique solution.

By similar considerations as in the previous subsection, we obtain a solution also in the congested case if condition \eqref{regular221} is not verified. 
To obtain a unique solution, we modify the priority rule described in \cite{garavello2016models} for a traffic flow. Here, the main goal is to use the whole capacity of the outgoing arc $i = 3$, which implies 
\begin{equation}
\label{eq:network_2to1_maximalitycond}
f_3 = \fmax_3 = a\sub{3} \rhomaxnet_3.
\end{equation}

\begin{figure}[htb]
	%		\centering
	\begin{center}
	\begin{minipage}{0.48\textwidth}
		\begin{tikzpicture} 
		\node[draw=none] at (-0.8, 3.5)() {$a)$} ;
		\draw[fill={rgb:black,.3;white,2}]	(0, 0) -- (3,0) -- (3, 1) -- (1.5, 2) -- (0, 2);
		\draw[->] (0,0) -- (4.0,0) node[right] {$f_1$};
		\draw[->] (0,0) -- (0,3.0) node[above] {$f_2$};
		\draw[-] (3,0) node[below] {$\fmax_1$} -- (3,1);
		\draw[dashed] (3,1) -- (3,2);
		\draw[-] (0,2) node[left] {$\fmax_2$} -- (1.5,2);
		\draw[dashed] (1.5,2) -- (3,2);
		\draw[-] (0,2) -- (1.5,2);
		\draw[-] (3.5, 0.66) -- (0.5, 2.66) node[above right] {$f_1 + f_2 = \fmax_3$} ;
		\draw[-] (0,0) --  (4,2) node[above] {$ f_2 = \frac{q}{1-q} f_1$};
		\draw[black,fill=black] (2.55, 1.275) circle (.5ex) node[above] {$P$};
		\node[draw=none] at (1, 1.3)(){$\Theta$};
		\end{tikzpicture}
	\end{minipage} 
%		\hspace{0.5cm}
	\begin{minipage}{0.48\textwidth}
		\begin{tikzpicture}
		\node[draw=none] at (-0.8, 3.5)() {$b)$} ;
		\draw[fill={rgb:black,.3;white,2}]	(0, 0) -- (3,0) -- (3, 1) -- (1.5, 2) -- (0, 2);
		\draw[->] (0,0) -- (4.0,0) node[right] {$f_1$};
		\draw[->] (0,0) -- (0,3.0) node[above] {$f_2$};
		\draw[-] (3,0) node[below] {$\fmax_1$} -- (3,1);
		\draw[dashed] (3,1) -- (3,2);
		\draw[-] (0,2) node[left] {$\fmax_2$} -- (1.5,2);
		\draw[dashed] (1.5,2) -- (3,2);
		\draw[-] (0,2) -- (1.5,2);
		\draw[-] (4, 0.33) -- (0.5, 2.66) node[above right] {$f_1 + f_2 = \fmax_3$} ;
		\draw[-] (0,0) --  (4,0.8);% node[above] {$ f_2 = \frac{q}{1-q} f_1$};
		\draw[black,fill=black] (3,1) circle (.5ex) node[below left] {$Q$};
		\draw[black,fill=black] (3.45,0.68) circle (.5ex) node[below] {$P$};
		\node[draw=none] at (1, 1.3)(){$\Theta$};
		\node[draw=none] at (4.2, 1.1)(){$ f_2 = \frac{q}{1-q} f_1$};
		\end{tikzpicture}
	\end{minipage}
	\end{center}
	\caption{Choice of the merging parameter $q$}
	\label{fig:merging_parameter}
\end{figure}
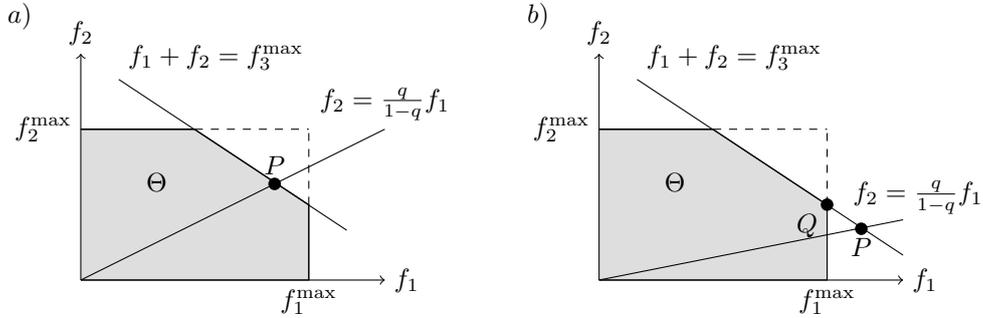

We set the merging parameter $q \in [0, 1]$ such that
\begin{equation}
\label{eq:network_2to1_fluxdist}
f_1 = q\fmax_3 \quad \text{ and } \quad f_2 = (1-q)\fmax_3.
\end{equation}
This leads to
\begin{equation}
\label{eq:network_2to1_ratiocond}
f_2 = \frac{1-q}{q} f_1,
\end{equation}
describing a ratio of the actual fluxes on the corresponding arcs.
The admissible region for the fluxes is 
\[ \Theta = \left\lbrace (f_1, f_2): 0 \leq f_1 \leq \fmax_1, 0 \leq f_2 \leq \fmax_2, 0 \leq f_1 + f_2 \leq \fmax_3 \right\rbrace \]
and shaded gray in Figure \ref{fig:merging_parameter}.
If condition \eqref{eq:network_2to1_maximalitycond} is not fulfilled by only considering the ratio \eqref{eq:network_2to1_ratiocond}, the parameter $q$ is adapted to obtain a unique solution. This is shown in Figure \ref{fig:merging_parameter}:
\begin{itemize}
	\item[a)]  The intersection point $P$ between the maximal outgoing flux (the line $f_1+f_2=\fmax_3$) and the priority ratio \eqref{eq:network_2to1_ratiocond} is inside the admissible set $\Theta$. In this case, we keep ${q \in (0,1)}$ fix and we have \eqref{eq:network_2to1_fluxdist}.
	\item[b)]  The intersection point $P$ is outside $\Theta$. We choose the closest point $Q$ inside $\Theta$ on the line $f_1 + f_2 = \fmax_3$, which guarantees maximal throughput, i.e., $f_1 = \fmax_1$. The merging parameter changes to $q=\fmax_1/\fmax_3.$ Then, the resulting fluxes are
	$f_1 = \fmax_1 $ and $ f_2 = \fmax_3 - \fmax_1.$
\end{itemize}

%\begin{remark}The same principle can be applied if arc $i$ does not provide enough mass to fulfill the maximality constraint $f_i = \fmax_i$. In this case, the domain $\Theta$ is adapted to the reduced actual mass flow $f_i = f_i^{\text{act}}$. If the domain $\Theta$ and the line $f_1 + f_2 = \fmax_3$ do not intersect, we are in the non-congested case described previously, and no merging rule is needed. 
%\end{remark}
We call $\Lambda_i$ the congested region (defined as in \eqref{laambda}) on arc $i=1,2$, $g_i(t)$ its interface and $q\subq{i} \in \{q, 1-q\}$ the corresponding merging parameters. The solution is then
\begin{equation}\label{sol221}
%\begin{multline}
\rho(x,t)=\\
\begin{cases}
\rhozeronetwork_i(x-a\sub{i} t), & (x,t)\in \Omega_i\setminus \Lambda_i \times (0,T],i=1,2\\
\rhomaxnet_i,
%\frac{a\sub{i}}{\bar a}\rho_i_0\left(g_i(\bar t_i)-a\sub{1}\left(t-\frac{x+g_i(\bar t_i)}{\bar a}\right)\right)
& (x,t)\in \Lambda_i \times (0,T],i=1,2\\
\rhomaxnet_3, & (x,t) \in \Omega_3 \times (0,T]:\begin{aligned}[t] & x\leq a\sub{3} t, \\ & \max\limits_{i=1,2}\Big\{g_i\big(t-\frac{x}{a\sub{3}}\big)\Big\}\neq 0 \end{aligned}\\
\sum\limits_{i=1,2}\frac{a\sub{i}}{a\sub{3}}\rhozeronetwork_i\left(-a\sub{i}\left(t-\frac{x}{a\sub{3}}\right)\right), & (x,t)\in \Omega_3 \times (0,T]:  \begin{aligned}[t] & x\leq a\sub{3} t ,\\  &\max\limits_{i=1,2}\Big\{g_i\big(t-\frac{x}{a\sub{3}}\big)\Big\} = 0 \end{aligned}\\
\rhozeronetwork_3(x-a\sub{3} t), & (x,t)\in \Omega_3 \times (0,T]: x>a\sub{3} t.
\end{cases}
%\end{multline}
\end{equation} 
If only one arc is congested, the solution holds true with $\Lambda_i = \emptyset$ for the non-congested arc $i$.
The shape of the congested region $\Lambda_i$ depends on the merging parameter $q\subq{i}$. It is described by the interface $g_i:[0, \infty) \rightarrow \R_{\leq 0}$, where $t$ is mapped to the solution of equation
\begin{equation}
-x+(t-t_0)\left(\frac{q\subq{i} \rhomaxnet_3 }{\rhomaxnet_i}a\sub{3} \right)-\frac{1}{\rhomaxnet_i}\int\limits_{x-a\sub{i} t}^{-a\sub{i} t_0}\rhozeronetwork_i(s)ds=0
\end{equation}
if this is negative, and to zero otherwise, analogously to the previous cases.

\begin{remark}
	The time $t_0$ is unique since congestion starts (independent on $q$) if the outgoing belt $3$ is not able to absorb all the incoming flux $f_1 + f_2$.
This is due to \eqref{eq:network_2to1_maximalitycond}, which ensures that congestion arises if the maximal capacity is reached. At the same time the evolution of the function $g_i(t)$ depends on the parameter $q$. Therefore, it is possible that $g_1(\hat t)=0$ for some $\hat t>t_0$ if $g_2(\hat t)>0$ (or vice-versa). 
This implies that for one arc $i \in \{1,2\}$, the function $g_i$ may start from zero while the other one starts from a negative value, i.e., at least one arc is congested, and the description of the $\Lambda_i$ is not completely separated.
\end{remark}
We now draw our attention the numerical treatment of the formerly stated problems.

\section{Numerical approximation}   \label{snum}

In this section, we present  
numerical experiments to illustrate and confirm our theoretical results for different network configurations.
The numerical scheme we propose is an adaptation of the scheme in \cite{towers2000convergence} to networks. 
This scheme has also been successfully applied to pedestrian networks in \cite{Camilli201793}. 

We discretize each arc $\Omega_i = (\intanf_i,\intend_i) \subset \R$ by
$X_i=(x_{i,0}=\intanf_i,x_{i,1},\dots, x_{i,\narc}=\intend_i)$
with constant discretization step $ \Delta x = {|x_{i,j} - x_{i,j-1}|}, ~ j = 1, \ldots, {\narc}$.
The spatial grid cells are defined as 
$\gridcell_{i,j} = (x_{i,j-1/2}, x_{i,j+1/2} ) \subset \mathbb{R},$
where the index $i$ refers to the corresponding arc.
We discretize the time set $[0,T]$ with $\{t^n=n\Delta t,\; n=0,1,\ldots,T/\Delta t\},$ where $\Delta t \in [0,T]$.
We define the piecewise constant approximation of the solution $\rho$ as
$\rhonum^n_{i,j}\approx\rho(x_{i,j},t^n)$
with $\rhonum^n$ constant in each grid cell $\gridcell_{i,j}$.
The scheme to update the approximation in each time step is for all $i \in E, j \in {2, \ldots, \narc-1}$ given by
\begin{equation}
\label{eq:scheme}
\left\{ \begin{array}{ll}
\rhonum^{n+1}_{i,j}= \rhonum^n_{i,j}-\frac{\Delta x}{\Delta t} \left( h\big(\rhonum^n_{i,j}, \rhonum^n_{i,j+1}\big) - h\big(\rhonum^n_{i,j-1},\rhonum^n_{i,j}\big) \right),\\
\rhonum^0_{i,j}=\rhozeronetwork(x_{i,j}). 
\end{array}
\right.
\end{equation}
For $j = 1$ and all $i \in E$, the update rule reads as 
\[\rhonum^{n+1}_{i,1}= \rhonum^n_{i,1}-\frac{\Delta x}{\Delta t} \left(h\big(\rhonum^n_{i,1}, \rhonum^n_{i,2}\big) - \hin_i  \right)\]
with $\hin_i$ describing an inflow condition at time $t^n$. Similarly, for $j = \narc$ and outflow condition $\hout_i$ at time $t^n$, the update rule is 
for all $i \in E$
\[\rhonum^{n+1}_{i,\narc}= \rhonum^n_{i,\narc}-\frac{\Delta x}{\Delta t} \left( \hout_i - h\big(\rhonum^n_{i,\narc-1},\rhonum^n_{i,\narc} \big) \right).\]
For the nodes $v \in V$ with $\delta_v^- = \emptyset$, we need an inflow condition which is set to zero, i.e.,
$\hin_{\iout} = 0 \text{ for } \iout \in \delta_v^+.$ 
For the nodes $v \in V$ with $\delta_v^+ = \emptyset$, we set the outflow as
$\hout_{\iin} = f_{\iin}(\rhonum^n_{\iin,\narcin}) \text{ for } \iin \in \delta_v^-.$
For all interior nodes $v \in V$ with $\delta_v^- \neq \emptyset$ and $\delta_v^+ \neq \emptyset$, we need to impose a junction rule at $v$
depending on the type of junction. 
The inflow and outflow conditions are defined as follows.
\begin{itemize}
    \item[a)] For an one-to-one junction, i.e.~$|\delta_v^-|=|\delta_v^+|=1$, we set
    \[
    \hout_{\iin} = \hin_{\iout} = h\big(\rhonum^n_{\iin,\narcin}, \rhonum^n_{\iout,1}\big) \text{ for } \iin \in \delta_v^-, \iout \in \delta_v^+. \]
    %This corresponds to a continuation of the numerical scheme through the junction point.
    \item[b)] For an one-to-two junction, i.e.~$1=|\delta_v^-|\neq|\delta_v^+|=2$, 
        we choose in the non-congested case \eqref{regular122}
    \[ \hout_{\iin} = f(\rhonum^n_{\iin,\narcin}) \text{ for }\iin \in \delta_v^-, \text{ and } \hin_{\iout} = \alpha_{\iout} f_{\iin}(\rhonum^n_{\iin,\narcin}) \text{ for }\iin \in \delta_v^-, \iout \in \delta_v^+  \]
    %and
    %\[ \hin_{\iout} = \alpha_{\iout} f_{\iin}(\rhonum^n_{\iin,\narcin}) \text{ for }\iin \in \delta_v^- \text{ and }\iout \in \delta_v^+.\]
    so that the flux is distributed according to the parameter $\alpha_{\iout}$.
    In the congested case, i.e., if \eqref{regular122} is violated, we consider the junction properties described in Section \ref{s122}.
    If we are in the case of a ``passive'' junction, we set the outgoing flux to
    \[ \hout_{\iin} = \sum_{\iout \in \delta_v^+} \hin_{\iout} \text{ for }\iin \in \delta_v^-, \]
    where the incoming fluxes on the outgoing arcs $\iout \in \delta_v^+= \{2,3\}$ are defined by
    \[ \hin_{2} = \min\Big\{ \fmax_{2}, \frac{\mu}{1-\mu} \fmax_3 \Big\}, \quad  \hin_{3} = \min\Big\{ \fmax_{3}, \frac{1-\mu}{\mu} \fmax_2 \Big\}. \] 
    If we consider an ``active'' junction instead, we set 
    \[ \hout_{\iin}  =\sum_{\iout \in \delta_v^+}\fmax_{\iout} \text{ for }\iin \in \delta_v^-, \text{ and } \hin_{\iout} = \fmax_{\iout} \text{ for } \iout \in \delta_v^+. \]
    %and the incoming fluxes to
    %\[ \hin_{\iout} = \fmax_{\iout} \text{ for } \iout \in \delta_v^+.\] 
    \item[c)] For a two-to-one junction, i.e.
    $2=|\delta_v^-|\neq|\delta_v^+|=1$, we choose in the non-congested case \eqref{regular221}
    \[ \hout_{\iin} = f_{\iin}(\rhonum^n_{\iin,\narcin}) \text{ for }\iin \in \delta_v^-, \text{ and } \hin_{\iout} = \sum_{\iin \in \delta_v^-} f_{\iin}(\rhonum^n_{\iin,\narcin}) \text{ for }\iout \in \delta_v^+. \]
    %and the incoming flux as sum of the outgoing fluxes:
    %\[ \hin_{\iout} = \sum_{\iin \in \delta_v^-} f_{\iin}(\rhonum^n_{\iin,\narcin}) \text{ for }\iout \in \delta_v^+.\] 
    In the congested case, we apply the merging parameter $q\subq{i}$ described in Section \ref{s221} and set
    \[ \hout_{\iin} = q\subq{\iin} \fmax_{\iout} \text{ for }\iin \in \delta_v^-, \iout \in \delta_v^+, \text { and } \hin_{\iout} = \fmax_{\iout} \text{ for }\iout \in \delta_v^+.\]
    %and
    %\[ \hin_{\iout} = \fmax_{\iout} \text{ for }\iout \in \delta_v^+.\] 
\end{itemize}

To define $h(\rhonum^n_{i,j},\rhonum^n_{i,j+1})$ in \eqref{eq:scheme}, we look for a function $h$ satisfying
\begin{eqnarray}
& h(0,0)=h(\rhomaxnet_i, \rhomaxnet_i)=0,        \label{flux0}\\
& m_-(\tilde{u})\le \frac{\partial}{\partial \tilde{u}} h(\tilde{u},u)\le 0\le \frac{\partial}{\partial u} h(\tilde{u},u)\le m_+(u),     \label{flux1}
\end{eqnarray}
with continuous function $m:\R \to\R$ and $m_-=\min(m,0)$, $m_+=\max (m,0)$. 
%In general, the maximal density $\rhomaxnet_i$ depends on the position within the network.
%However, we will state a flux function, for which equality \eqref{flux0} holds true everywhere. 

\begin{figure}[!t]
    \begin{center}
        \begin{tikzpicture}    
        \put(0,0){\node[draw=none]{\includegraphics[width=0.5\textwidth]{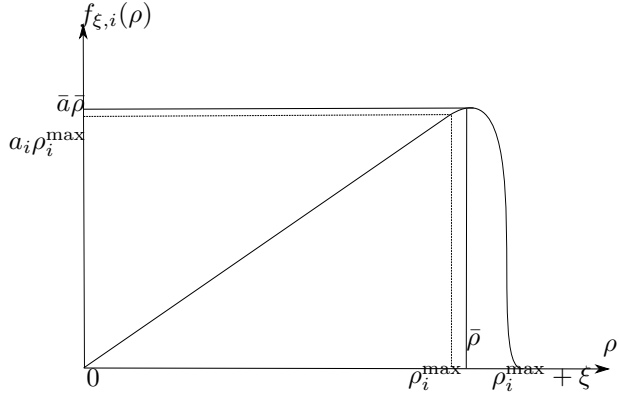}};}
        \put(0,0){\node[draw=none] at (-3.3,-2.3) {$0$};}
        \put(0,0){\node[draw=none] at (3.5, -1.9) {$\rho$};}
        \put(0,0){\node[draw=none] at (2.6, -2.3) {$\rhomaxnet_i + \smoothpar$};}
        \put(0,0){\node[draw=none] at (1.2, -2.3) {$\rhomaxnet_i$};}
        \put(0,0){\node[draw=none] at (-3, 2.5) {$f_{\smoothpar, i}(\rho)$};}
        \put(0,0){\node[draw=none] at (-3.6, 1.3) {$\bar{a} \bar{\rho}$};}
        \put(0,0){\node[draw=none] at (-3.9, 0.8) {$a\sub{i} \rhomaxnet_i$};}
        \put(0,0){\node[draw=none] at (1.7, -1.8) {$\bar{\rho}$};}
        \end{tikzpicture}
        \caption{Regularized flux function $f_{\smoothpar, i}$}
        \label{prob}
    \end{center}
\end{figure}

Conditions \eqref{flux0} and \eqref{flux1} enable to use results from \cite{Camilli201793,towers2000convergence}. 
Since $f_i$ is discontinuous in $\rho$, we need a suitable regularization.
We define a \textit{Friedrichs mollifier} $\varphi \in C_0^{\infty}(\R)$ with compact support in $[-1,1]$ such that 
\[\varphi(-y)=\varphi(y), \quad \int_\R \varphi(y)dy=1.\]
In our case, we use the mollifier ${\varphi(y):=\max(0,1-|y|)}$ and define 
$\varphi_\smoothpar(y):=\frac{2}{\smoothpar}\varphi(\tfrac{2y}{\smoothpar})$
for a small parameter $\smoothpar>0$, which implies that $\varphi_\smoothpar(y)$ has compact support in $[-\frac{\smoothpar}{2},\frac{\smoothpar}{2}]$.
For each arc $i \in E$, we introduce the following smooth regularization of the flux function \eqref{eq:network_flux_arci} 
\begin{equation}\label{regf}
f_{\smoothpar,i}(\rho):=a\sub{i} \rho\left(1-\int_{\rhomaxnet_i}^{\rho}\varphi_{\smoothpar}\left(y-\rhomaxnet_i-\frac{\smoothpar}{2}\right)dy\right),
\end{equation}
see Figure \ref{prob}. The function coincides with the original flux function \eqref{eq:network_flux_arci} in $x\in[0,\rhomaxnet_i]$ but there is a continuously differentiable connection to the value ${f_{\smoothpar,i}(\rhomaxnet_i+\smoothpar)=0}$, i.e., the regularized flux function $f_{\smoothpar,i}$ itself is continuously differentiable. We note that
$(f_{\smoothpar,i}(\rhomaxnet_i+\smoothpar))'= 0$ and $(f_{\smoothpar,i}(\bar \rho))'=\bar{a},$
i.e.~the transport velocity inside the congested area $\Lambda$. Moreover, the derivative is bounded by
$|(f_{\smoothpar,i})'(\rho)|\leq \frac{2}{\smoothpar}$
for small $\smoothpar$.
In the limit $\smoothpar\rightarrow 0^+$, we recover the original discontinuous flux function \eqref{eq:network_flux_arci}.

Note that the second equality of condition \eqref{flux0} translates to 
\begin{equation}
\label{eq:network_fluxreg_max}
h(\rhomaxnet_i + \smoothpar, \rhomaxnet_i + \smoothpar) = 0
\end{equation}
in the regularized case.

We choose the numerical flux function $h$ as \emph{Godunov flux}
\begin{equation}\label{eq:network_Godunov}
h(\rhonum^n_{i,j},\rhonum^n_{i,j+1})=\left\{
\begin{array}{ll}
\min\limits_{z \in [\rhonum^n_{i,j},\rhonum^n_{i,j+1}]} f_{\smoothpar,i}(z), & \hbox{if $\rhonum^n_{i,j}\le \rhonum^n_{i,j+1}$} \\[5pt]
\max\limits_{z \in [\rhonum^n_{i,j+1},\rhonum^n_{i,j}]} f_{\smoothpar,i}(z), & \hbox{if $\rhonum^n_{i,j}\ge \rhonum^n_{i,j+1}$}.
\end{array}
\right.
\end{equation}
and condition \eqref{flux1} is then satisfied with
$m(\rho)=\big(f_{\smoothpar,i}\big)'(\rho).$

The scheme \eqref{eq:scheme} is stable, if the following CFL condition 
\begin{equation}
\Delta t\leq \frac{\Delta x} {\max\limits_{v \in V}|\delta_v^-| \cdot \|m\|_{L^\infty(0,\rhomaxnet + \smoothpar)}} \label{CFL}
\end{equation}
is fulfilled, cf.~\cite{towers2000convergence}. 
%Moreover, therein it is shown, that the computed solution remains within $[0, \rhomaxnet_i]$, which means that the solution is physically reasonable.\\
From inequality \eqref{CFL}, we can establish a relation between the regularization parameter $\smoothpar$ and the discretization steps $\Delta t$ and $\Delta x$. In particular, if $\max|\delta_v^-|\leq 2$ for a fixed $v \in V$,
the CFL condition reduces to
\begin{equation}
\label{eq:Network_CFL}
\frac{\Delta t}{\Delta x}\leq \frac{\smoothpar}{4}.
\end{equation}
Knowing that $\smoothpar$ is supposed to be small, this is a quite restrictive condition. However, in the next section we will see that the 
choice of very small parameters $\smoothpar$ does not improve the solution significantly. 

%%%%%%%%%%%%%%%%%%%%%%%%%%%

\section{Tests}    \label{stests}
In this section, we present numerical results for different network settings. 
%On the one hand, we give a numerical validation to the analytical solutions discussed before, and on the other hand, we show that the proposed model is able to capture most of the non-trivial effects which are present in this kind of complex systems. 
Throughout this section, we compute the numerical solution  based on scheme \eqref{eq:scheme} with Godunov flux \eqref{eq:network_Godunov}.
If not stated otherwise, the space step size is $\Delta x=5\cdot 10^{-3}$, the time step size $\Delta t=10^{-5}$
and the smoothing parameter is fixed to $\smoothpar=10^{-2}.$ For simplicity, we choose $\rhomaxnet=1$ in all experiments.

\begin{figure}[t]
	\centering
%	\subfloat[$t=0$\label{fig:network_test1_t00}]
%	{
	\begin{minipage}{4.1cm}
		\includegraphics[width=1\textwidth]{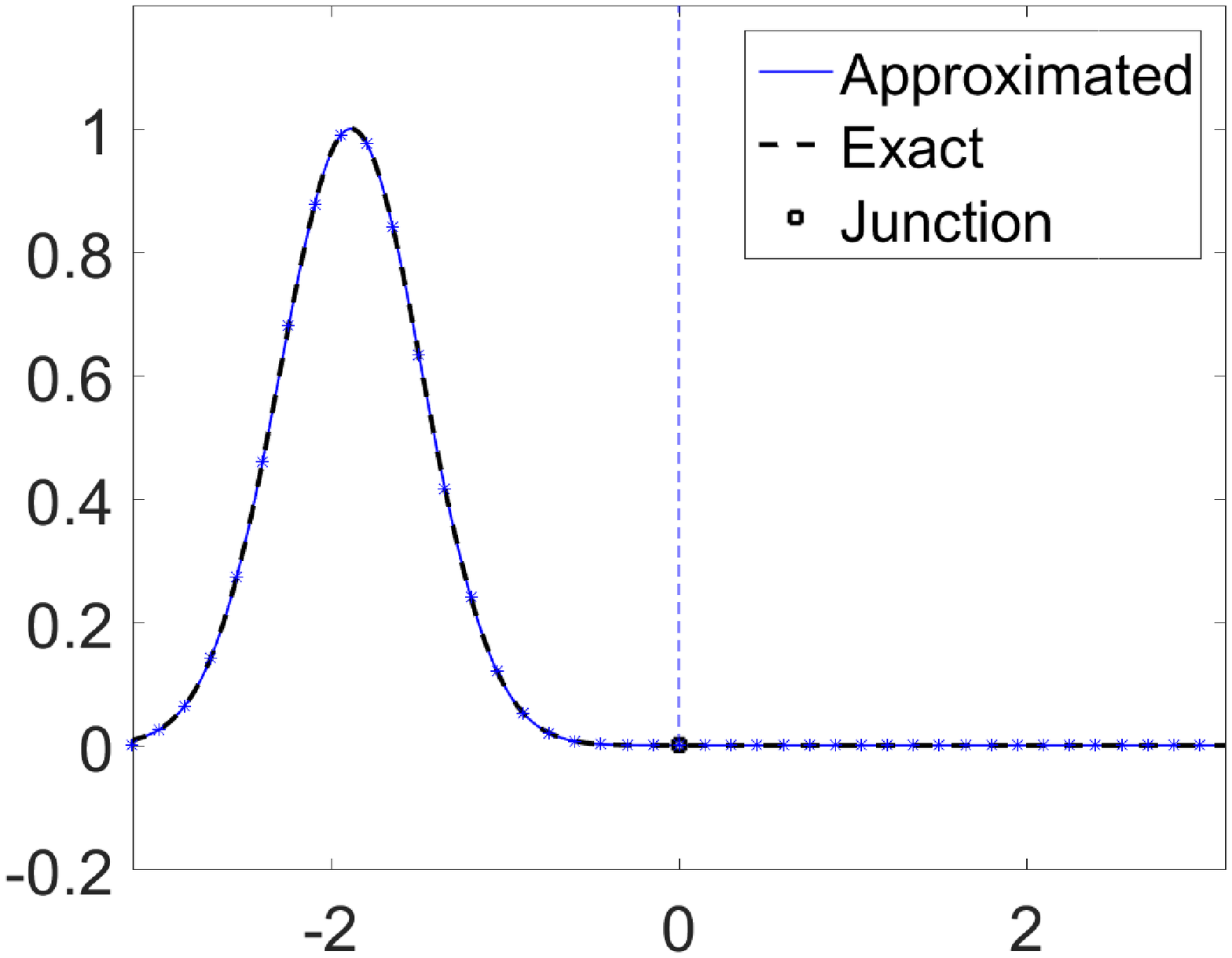}
		\centering
		\scriptsize $t = 0$
	\end{minipage}
%	}
%	\subfloat[$t=1$]
%	{
	\begin{minipage}{4.1cm}
		\includegraphics[width=1\textwidth]{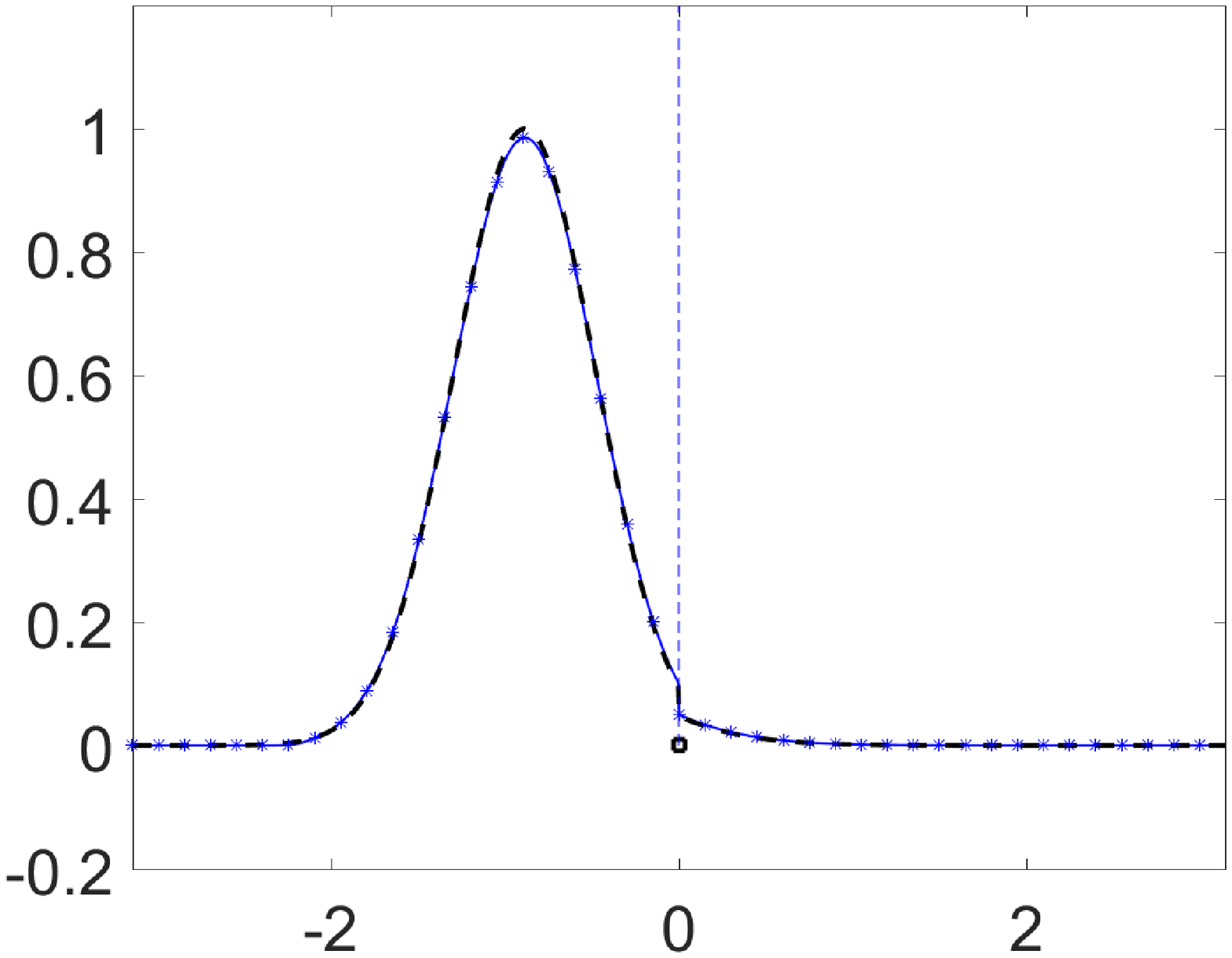}
		\centering
		\scriptsize $t = 1$
	\end{minipage}
%	} \\
%	\subfloat[$t=1.4$]
%	{
	\begin{minipage}{4.1cm}
		\includegraphics[width=1\textwidth]{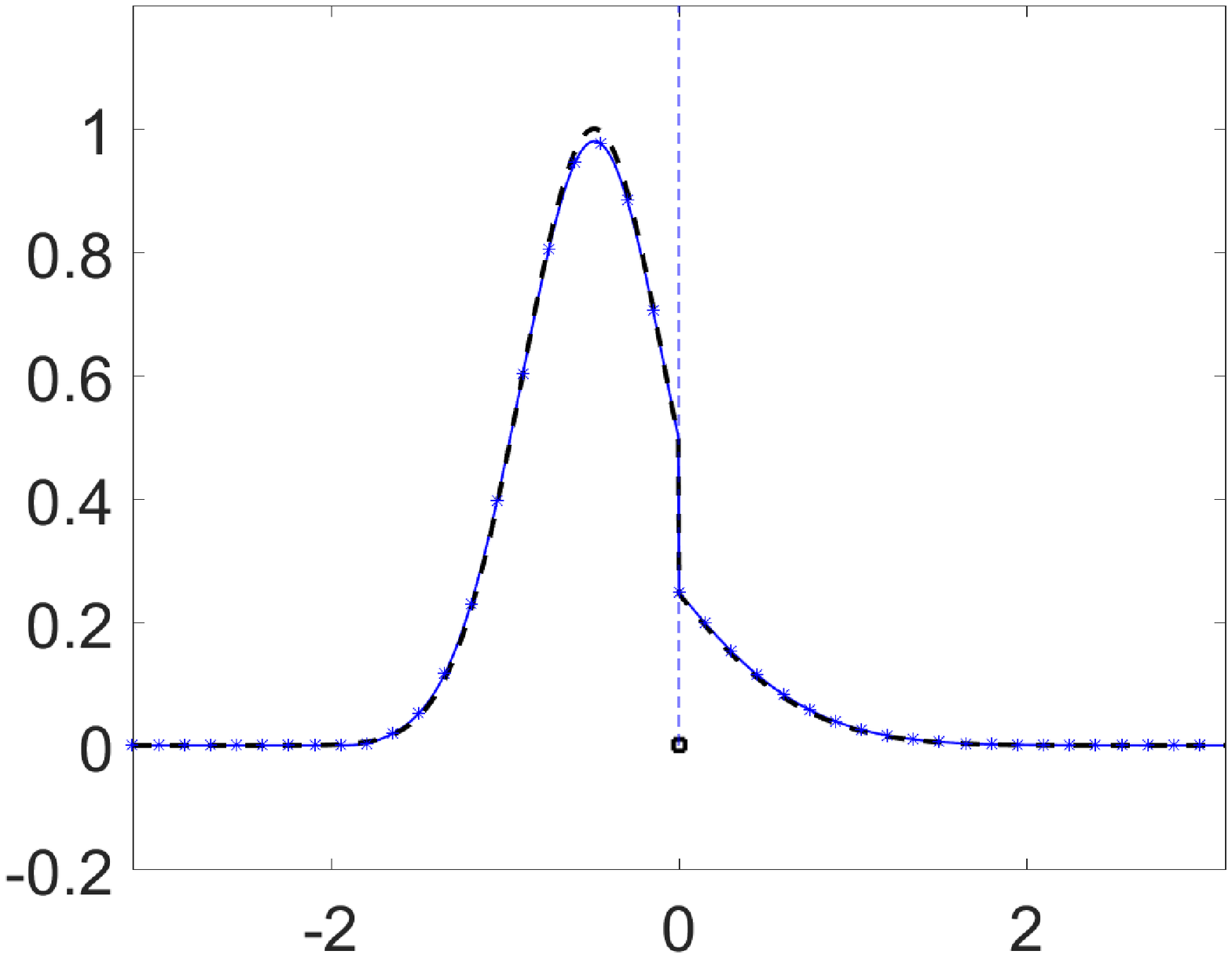}
		\centering
		\scriptsize $t = 1.4$
	\end{minipage}\\
%	}
%	\subfloat[$t=1.8$]
%	{
	\begin{minipage}{4.1cm}
		\includegraphics[width=1\textwidth]{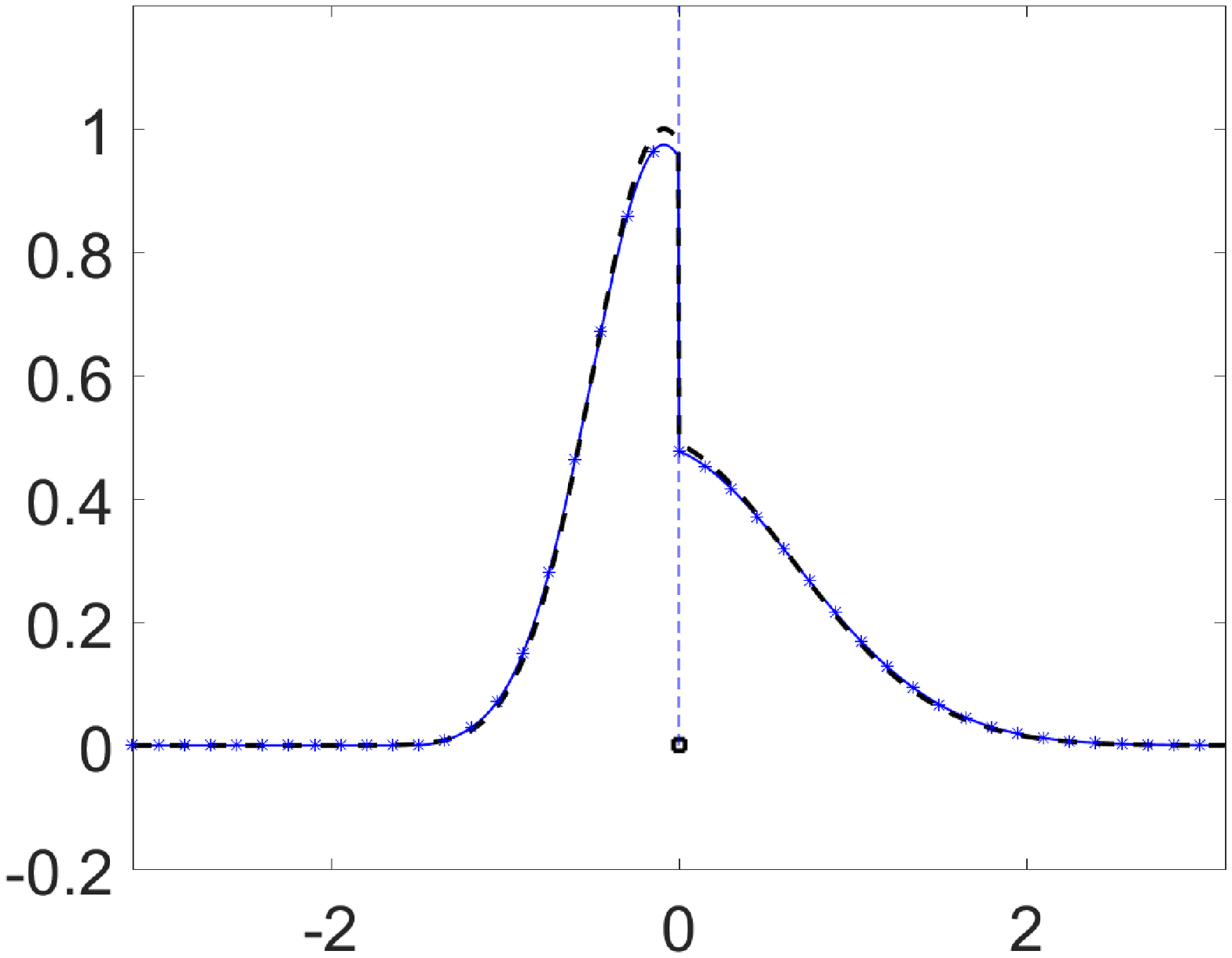}
		\centering
		\scriptsize $t = 1.8$
	\end{minipage}
%	} \\
%	\subfloat[$t=2.2$]
%	{
	\begin{minipage}{4.1cm}
		\includegraphics[width=1\textwidth]{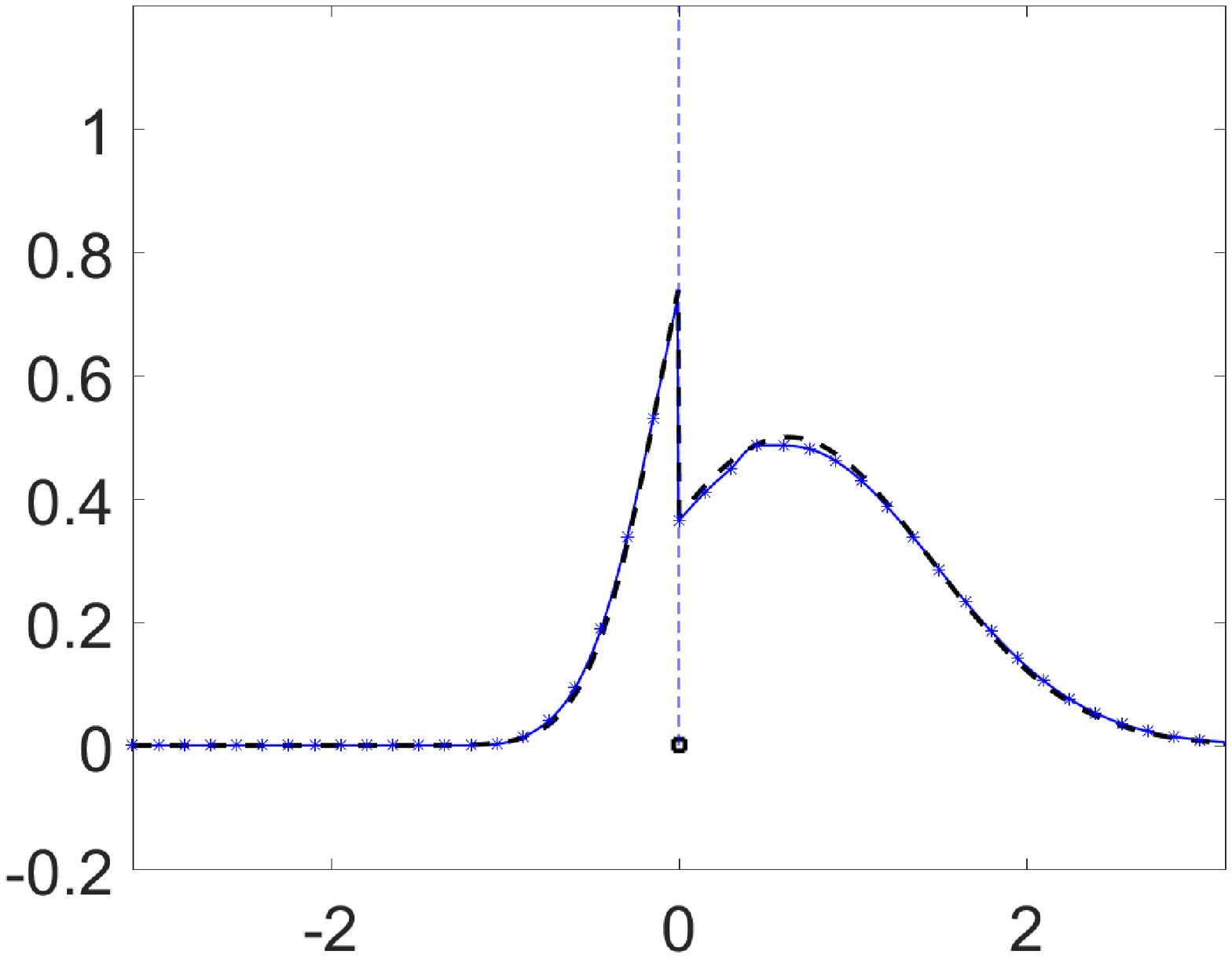}
		\centering
		\scriptsize $t = 2.2$
	\end{minipage}
%	}
%	\subfloat[$t=2.6$\label{fig:network_test1_t26}]
%	{
	\begin{minipage}{4.1cm}
		\includegraphics[width=1\textwidth]{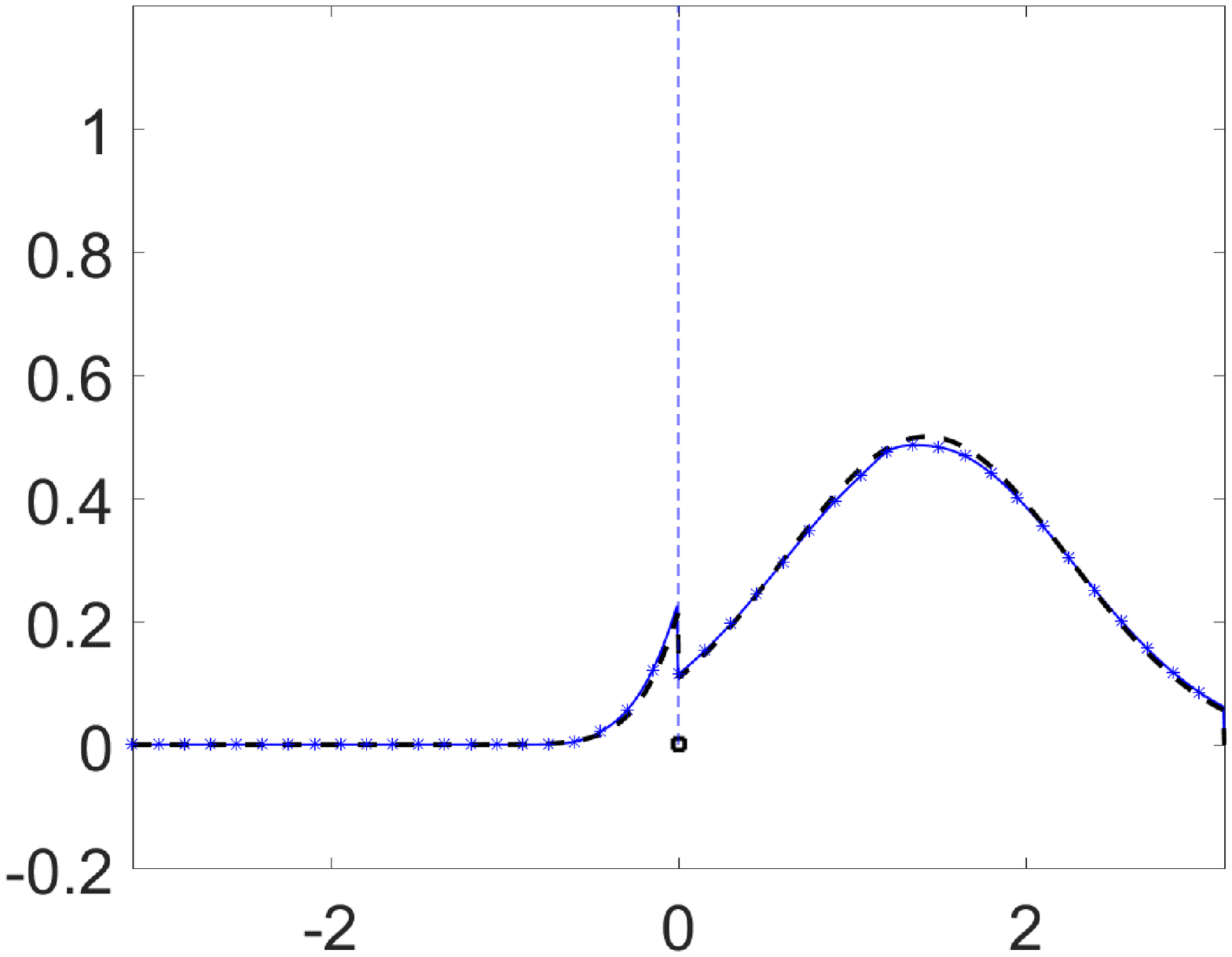}
		\centering
		\scriptsize $t = 2.6$
	\end{minipage}
%	}
	\caption{Test 1: non-congested case with $a\sub{1}=1$ and $a\sub{2}=2$}
	\label{fig:network_test1}
\end{figure}

\subsection{One-to-one junction} 
%\begin{table}[!b]
	%\centering
	%\begin{tabular}{lcc}
%%		\toprule
		%\multicolumn{2}{c}{Parameter} & Value \\ \hline
		%Velocity on arc $1$& $a\sub{1}$ & 1\\
		%Velocity on arc $2$& $a\sub{2}$ & 2\\
		%Maximal density on arc $1$& $\rhomaxnet_1$ & 1\\
		%Maximal density on arc $2$& $\rhomaxnet_2$ & 1\\
		%Space step size & $\Delta x$ & $5\cdot 10^{-3}$\\
		%Time step size & $\Delta t$ & $10^{-5}$\\
		%Smoothing parameter &$\smoothpar $& $10^{-2}$
%%		\bottomrule
	%\end{tabular}
	%\caption{Parameters for Test 1 and 2 on the one-to-one junction}
	%\label{tab:network_test1_parameters}
%\end{table}
First, we study the situation described in Section \ref{s121}. The linear network is given by $\Omega_1=(-\pi,0)$ and $\Omega_2=(0,\pi)$, i.e.~the coupling is at $x=0$. We fix the initial solution $\rhozeronetwork$ on $\Omega_1 \cup \Omega_2$ as 
\begin{equation}
\label{eq:network_test_rhozero}
\rhozeronetwork(x)=\exp\left({-3\Big(x+\frac{3}{5}\pi\Big)^2}\right).
\end{equation}

\subsubsection*{Test 1: free-flow case}

This is the non-congested case, i.e.~condition \eqref{regular} is holds true and the analytical solution is simply \eqref{sol:reg}.
The results are displayed in Figure \ref{fig:network_test1}. It shows the evolution of the initial density at different time steps
for the different velocities $a\sub{1}=1$ and $a\sub{2}=2$.
We see that the analytical and the numerical solution match very well. A discontinuity appears in the solution at $x = 0$. There, the doubling of the velocity  has the effect of ``spreading'' the initial solution. Due to the mass conservation, the local density behind the junction point is halved. This can be seen by comparing the maximal arising density $\rho = 1$ in front of the junction and the maximal arising density behind the junction which is $\rho = 0.5$.

\subsubsection*{Test 2: congested case}

\begin{figure}
	\centering
%	\subfloat[$t=0$\label{fig:network_test2_t00}]
%	{
	\begin{minipage}{4.1cm}
		\includegraphics[width=1\textwidth]{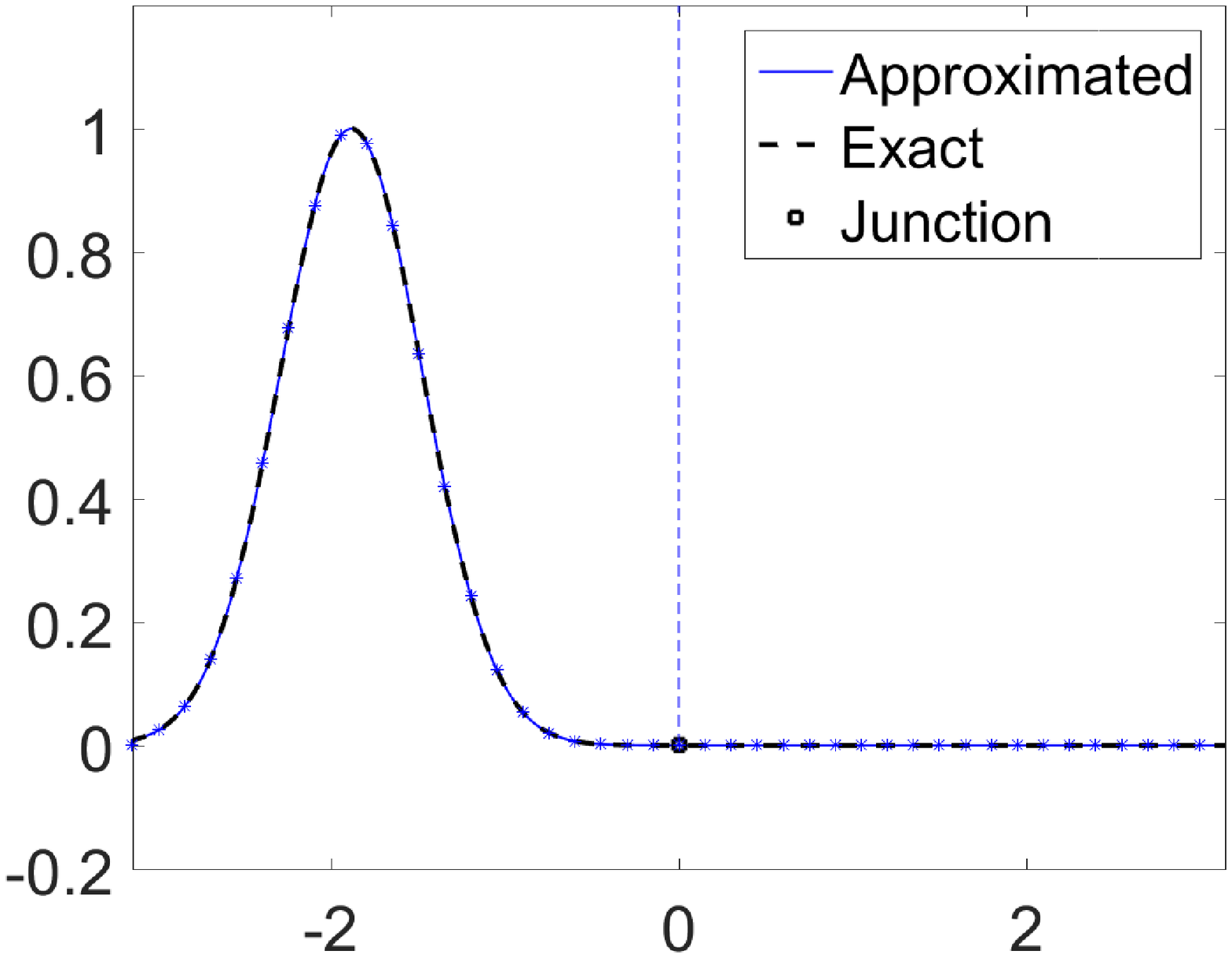}
		\centering
		\scriptsize $t=0$
	\end{minipage}
%	}
%	\subfloat[$t=0.5$\label{fig:network_test2_t05}]
%	{
	\begin{minipage}{4.1cm}
		\includegraphics[width=1\textwidth]{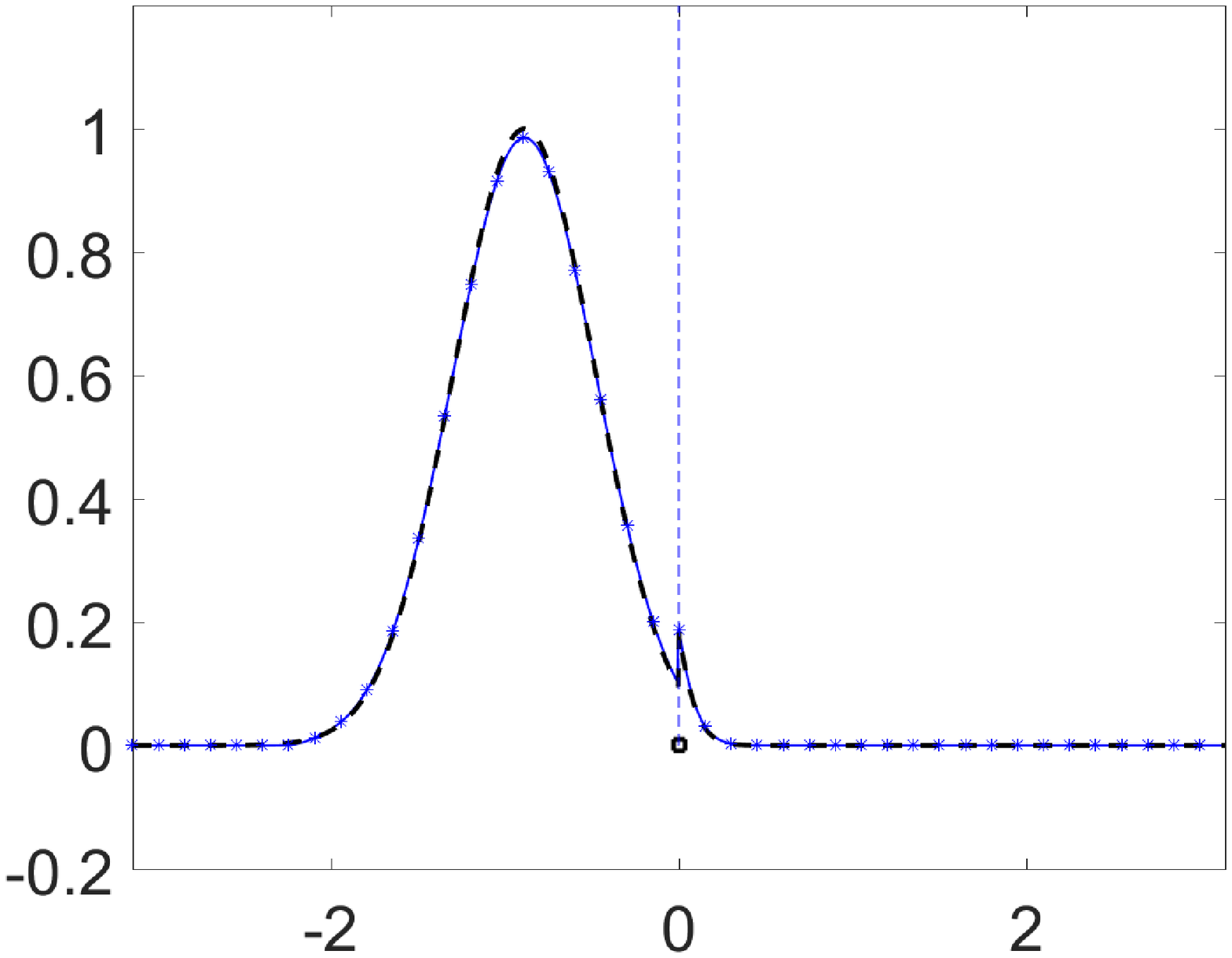}
				\centering
		\scriptsize $t=0.5$
	\end{minipage}
%	}
%	\subfloat[$t=0.8$\label{fig:network_test2_t08}]
%	{
	\begin{minipage}{4.1cm}
		\includegraphics[width=1\textwidth]{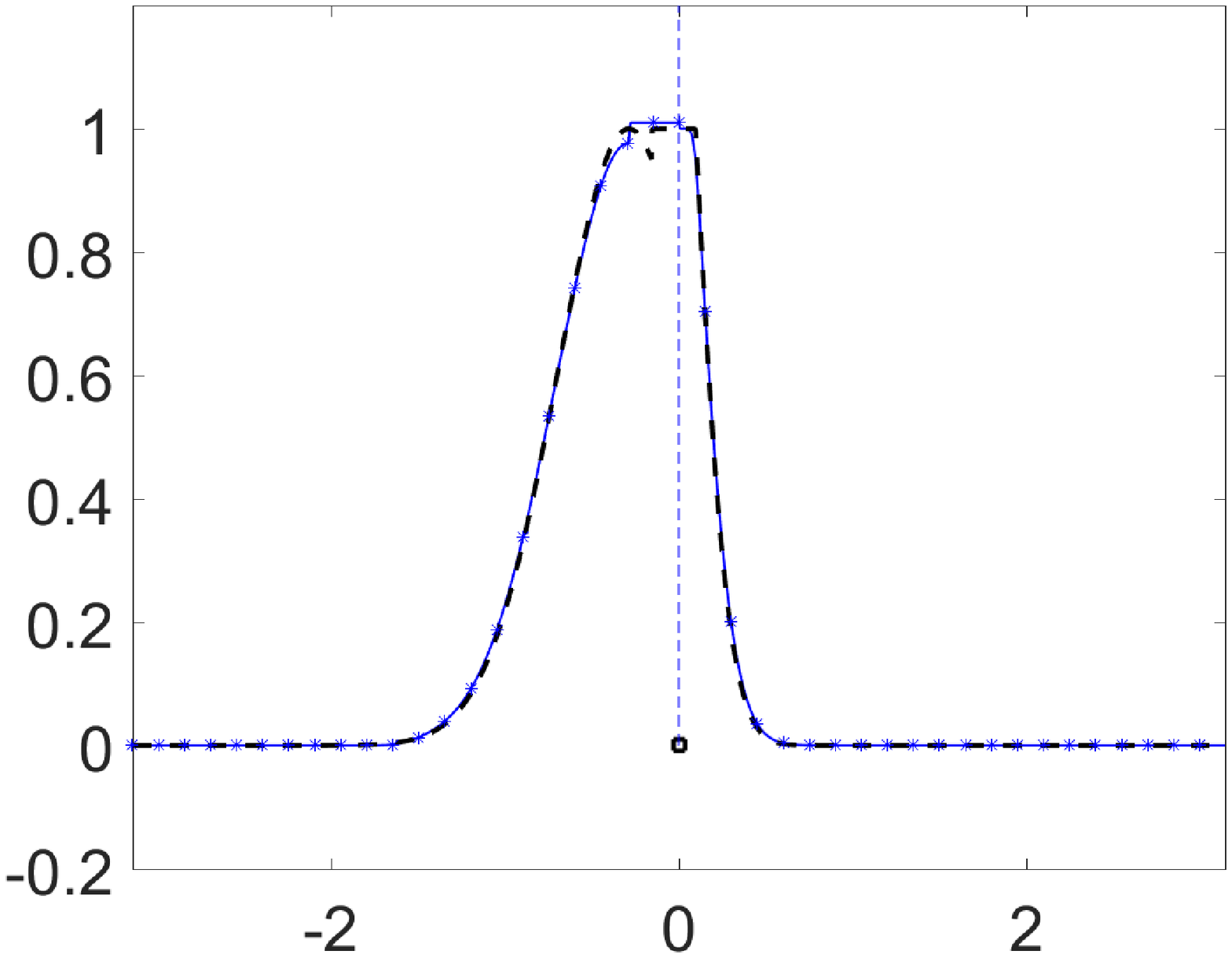}
				\centering
		\scriptsize $t=0.8$
	\end{minipage}
%	}
%	\subfloat[$t=1.1$\label{fig:network_test2_t11}]
%	{
	\begin{minipage}{4.1cm}
		\includegraphics[width=1\textwidth]{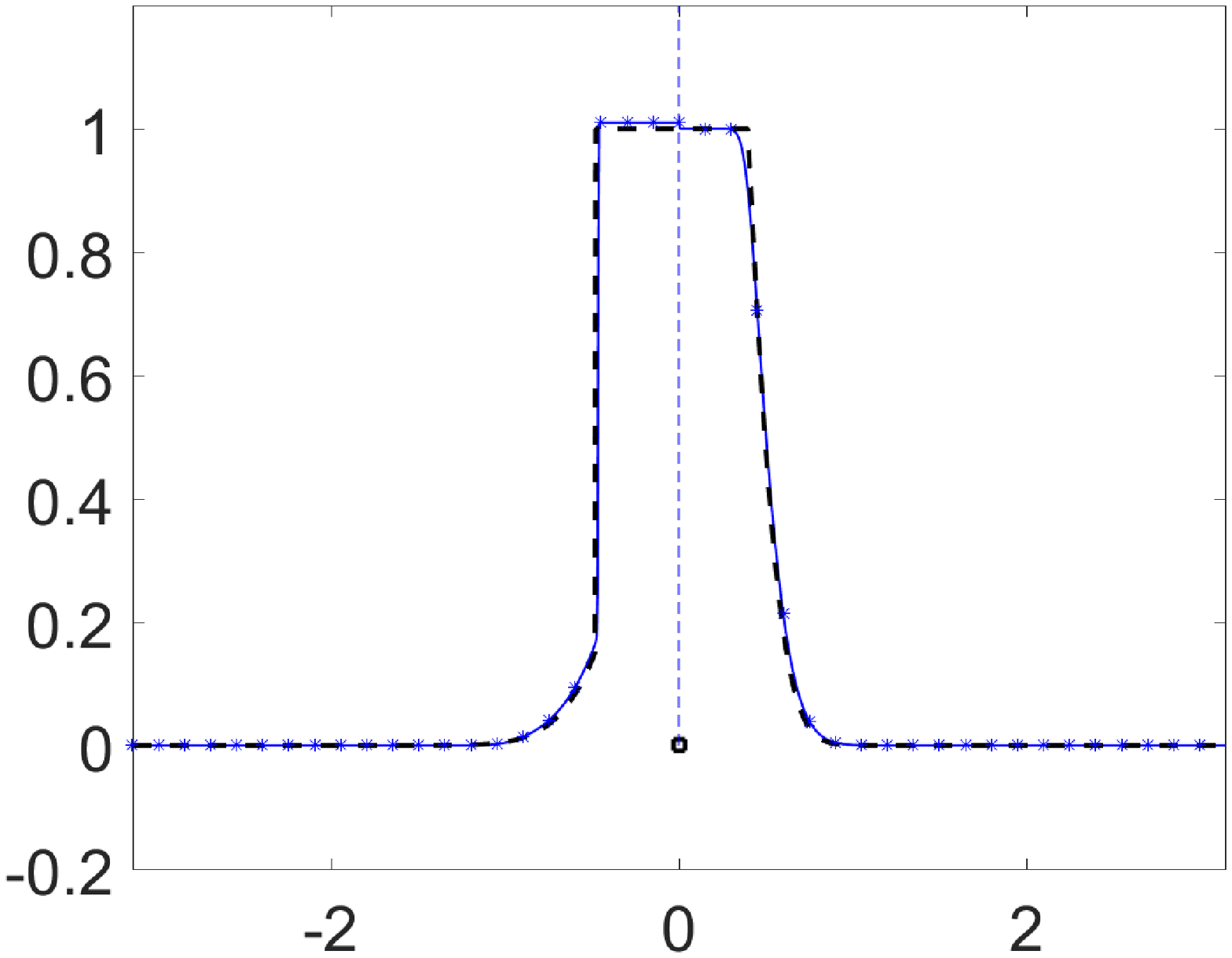}
				\centering
		\scriptsize $t=1.1$
	\end{minipage}
%	}\\
%	\subfloat[$t=1.4$\label{fig:network_test2_t14}]
%	{
	\begin{minipage}{4.1cm}
		\includegraphics[width=1\textwidth]{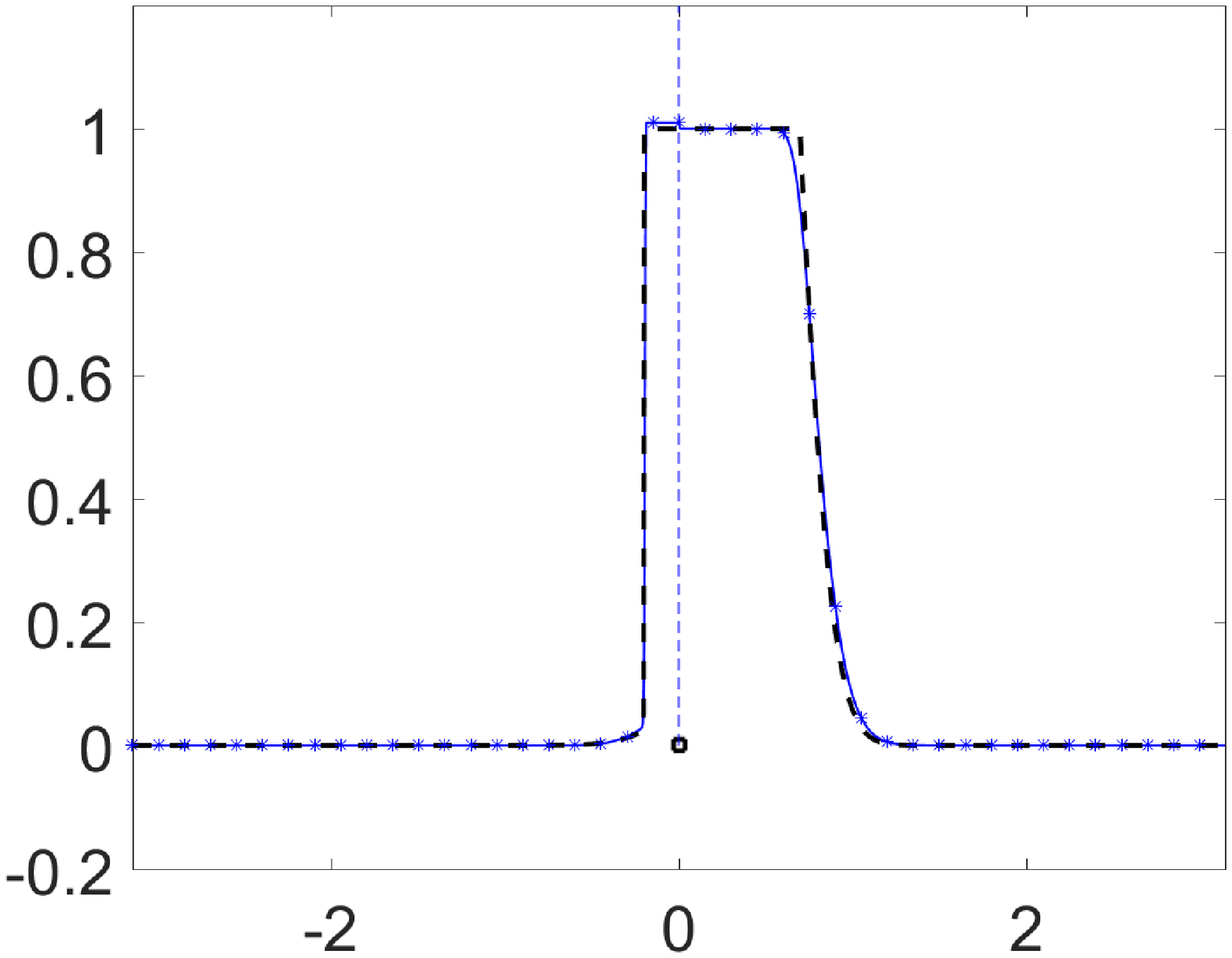}
				\centering
		\scriptsize $t=1.4$
	\end{minipage}
%	}
%	\subfloat[$t=1.7$\label{fig:network_test2_t17}]
%	{
	\begin{minipage}{4.1cm}
		\includegraphics[width=1\textwidth]{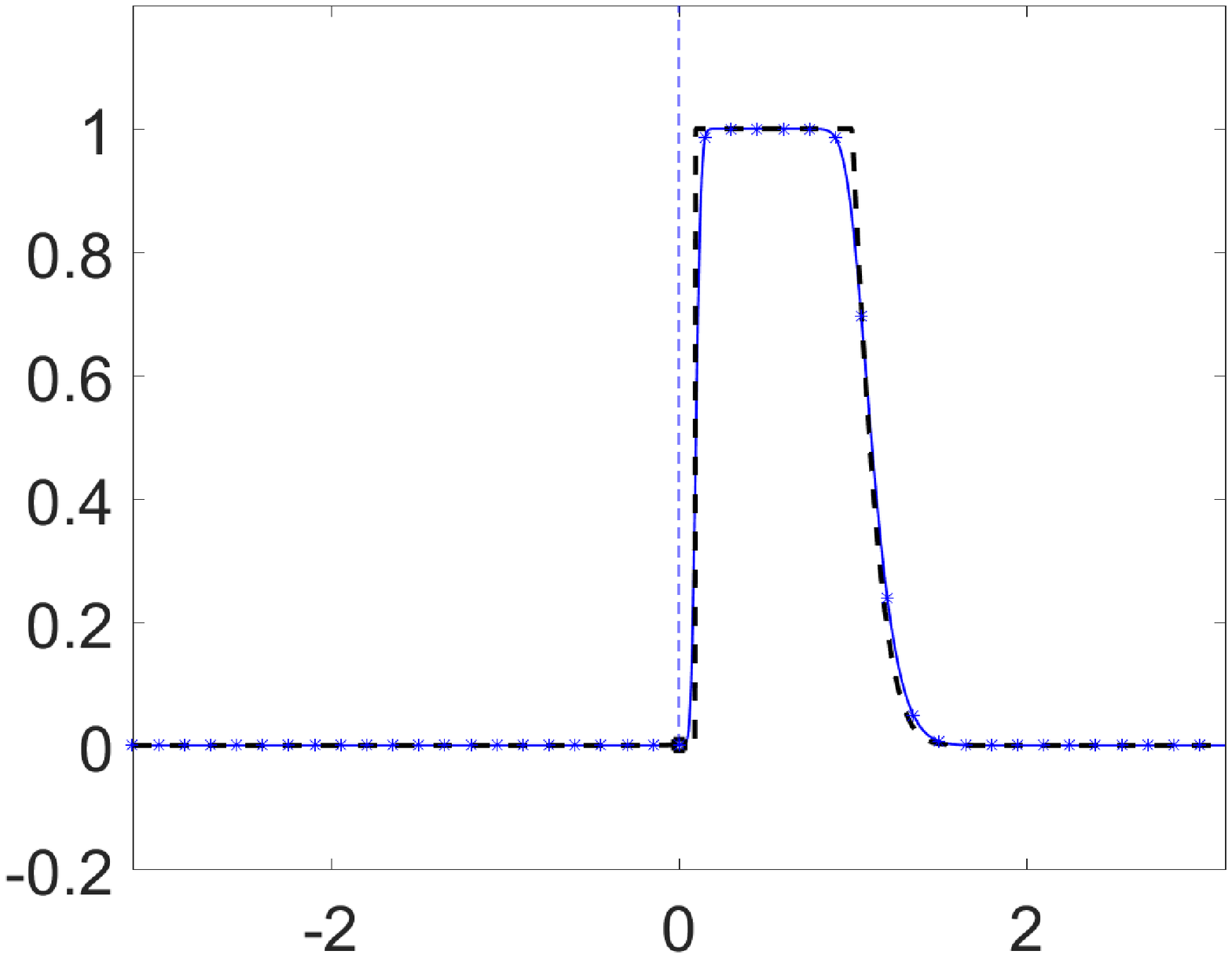}
				\centering
		\scriptsize $t=1.7$
	\end{minipage}
%	}\\
	\caption{Test 2: congested case with $a\sub{1}=2$ and $a\sub{2}=1$}
	\label{fig:network_test2}
\end{figure}

For the second test, we only vary the velocity of the arcs and we set $a\sub{1}=2$ and $a\sub{2}=1$. 
%All other numerical parameters remain the same as in the previous case in Table \ref{tab:network_test1_parameters}. 
Then, the condition \eqref{regular} is not true anymore and we are in the congested case, where the analytical solution is given by \eqref{sol121}.
Figure \ref{fig:network_test2} shows the comparison between the analytical and numerical solution.
After the time $t = 1.7$, the evolution continues on arc $2$ as linear transport.
Obviously, the discontinuity backward wave, modeled by the function $g(t)$ in \eqref{sol121}, is correctly tracked by the numerical solution. It should be noticed that the density located in the congested region $\Lambda$ in $\Omega_1$ moves with the velocity $\bar{a} = a\sub{2}$, cf.~\eqref{eq:network_1to1_bara}. 

\begin{figure}[htb]
	\begin{center}
%		\subfloat[Numerical solution\label{Fig:network_test2_num_spacetime}]{
		\begin{minipage}{5cm}
			\includegraphics[width=1\textwidth]{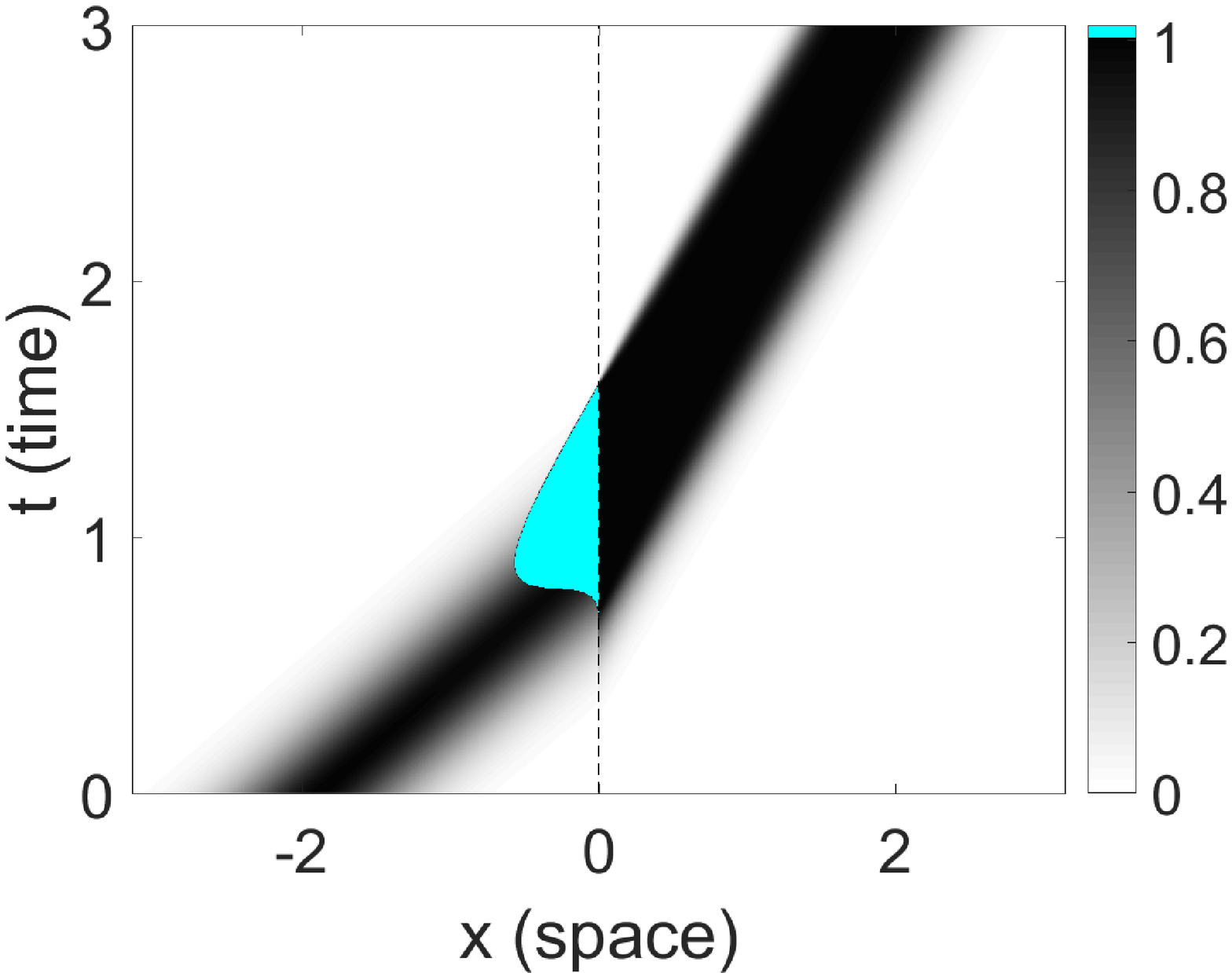}
			\centering
			\scriptsize 
			\textit{Numerical solution}
		\end{minipage}
%		}
%		\subfloat[Analytical solution\label{Fig:network_test2_anal_spacetime}]{
		\begin{minipage}{5cm}
			\includegraphics[width=1\textwidth]{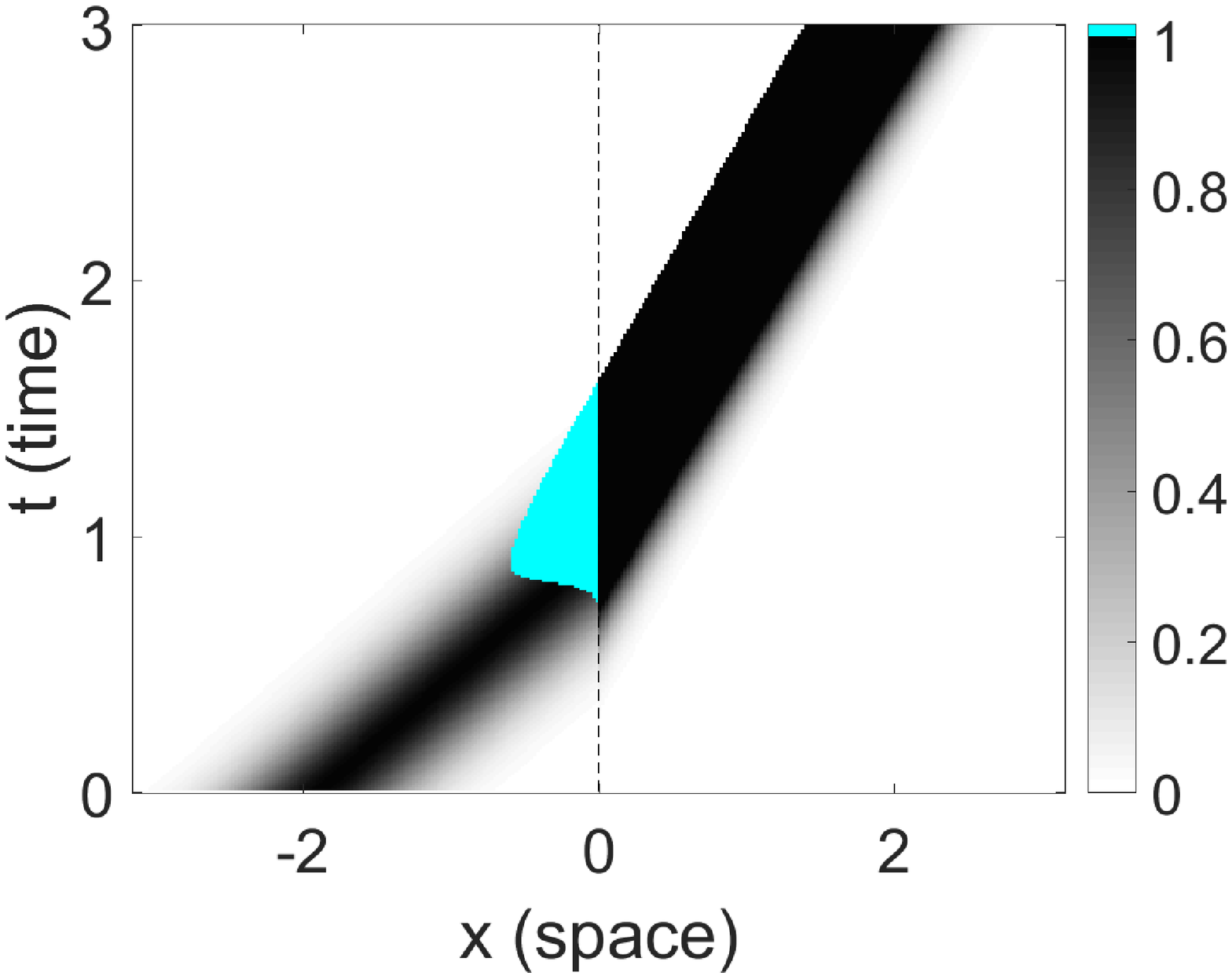}
			\centering
			\scriptsize \textit{Analytical solution}
		\end{minipage}
%		}
		\caption{Test 2: space-time diagram for the congested case}
		\label{t2f}
	\end{center}
\end{figure}

Figure \ref{t2f} shows the space-time diagram for the numerical and analytical solution.
Note that the highlighted congested region $\Lambda$ in the middle is correctly tracked by the numerical scheme. 
However, we observe the diffusive effect of the Godunov scheme. The latter effect could be reduced by the use of 
other flux approximations (as e.g.~proposed in \cite{fjordholm2012arbitrarily}) but not avoided completely. 

For the same setting, the influence of the discretization step sizes is shown in Table \ref{tab:network_analyze_delta} (left). Since the analytical solution is known, we can evaluate the $L^2$-error to the numerical solution. 
According to the space step size, the time step size is adapted to satisfy the CFL condition \eqref{CFL}.
As expected, the error tends to zero with decreasing step sizes. 

To study the influence of the smoothing parameter $\smoothpar$, the fixed discretization is chosen such that the CFL condition \eqref{eq:Network_CFL} is fulfilled for the smallest value of $\smoothpar$. This leads to a time step ${\Delta t = 2\cdot 10^{-6}}$ with a space step size $\Delta x = 5\cdot 10^{-3}$. The result is shown in Table \ref{tab:network_analyze_delta} (right). The error turns out to be only slightly smaller for the smallest value of $\smoothpar$. 
\begin{table}[tb] 
	\centering
	\begin{tabular}{ccc}
		%		\toprule
		$\Delta x$ & $\Delta t$ & error\\ \hline \\[-0.4cm]
		0.1\phantom{00}  & $2 \cdot 10^{-4}$ & 0.0842 \\
		0.05\phantom{0} &$\phantom{2 \cdot }10^{-4}$ & 0.0381 \\
		0.02\phantom{0} &$5 \cdot 10^{-5}$ & 0.0184 \\
		0.01\phantom{0} &$ 2 \cdot 10^{-5}$ & 0.0073 \\
		0.005 &$\phantom{2 \cdot }10^{-5}$ & 0.0057 
		%		\bottomrule
	\end{tabular}
\hspace{1cm}
	\begin{tabular}{cc}
		%		\toprule
		$\smoothpar$ & error\\ \hline \\[-0.4cm]
		$5 \cdot 10^{-2}$ & 0.0051 \\
		$2 \cdot 10^{-2}$  & 0.0042 \\
		$\phantom{1 \cdot} 10^{-2}$& 0.0039 \\
		$5 \cdot 10^{-3}$& 0.0037 \\
		$2 \cdot 10^{-3}$ &0.0035 
		%		\bottomrule
	\end{tabular}
%	\scriptsize \textit{Decreasing smoothing parameter}
	\caption{Decreasing step sizes (left), decreasing smoothing parameter $\smoothpar$ (right)}
	\label{tab:network_analyze_delta}
\end{table}

\begin{figure}[t!]
	\centering
	\begin{minipage}[b]{0.13\textwidth}
		$t = 0:$
		\vfill
	\end{minipage}
	\begin{minipage}[]{0.75\textwidth}
		\includegraphics[trim = 50 140 50 130, clip,
		width=1\textwidth]{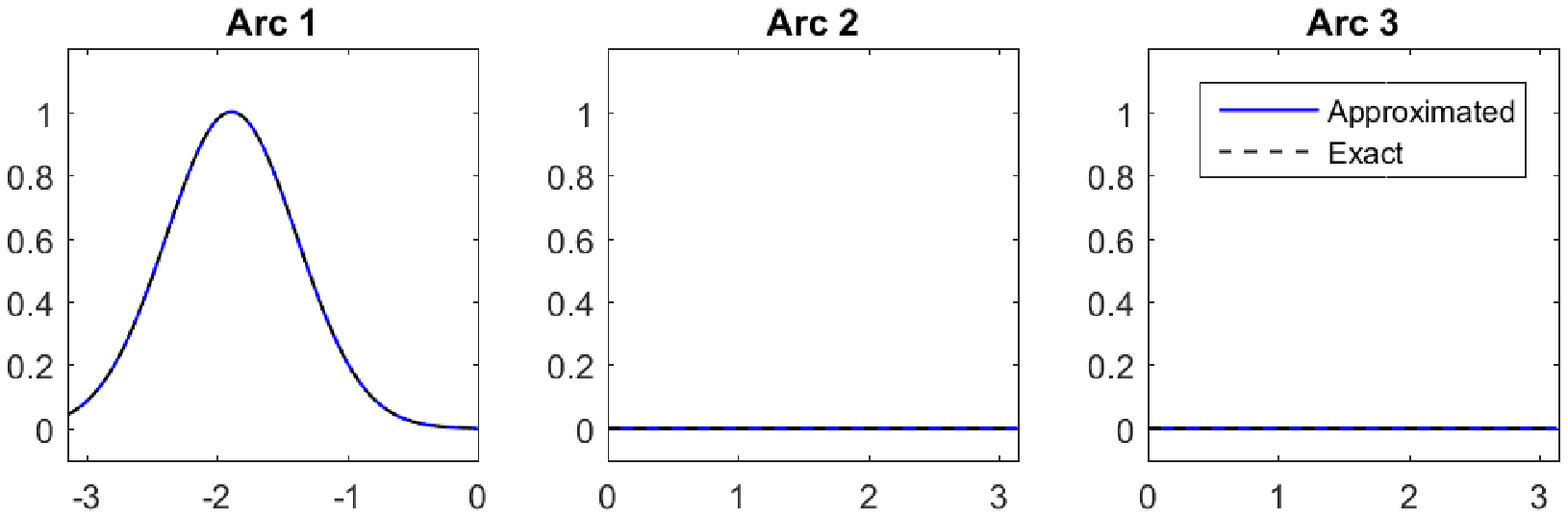}
	\end{minipage}\\
	%	\begin{minipage}[b]{0.13\textwidth}
	%		t = 0.375:
	%		\vfill
	%	\end{minipage}
	%	\begin{minipage}[]{0.77\textwidth}
	%		\includegraphics[trim = 50 140 50 130, clip, width=1\textwidth]{Pictures/Network/112_k0375_passive.eps}
	%	\end{minipage}\\
	\begin{minipage}[b]{0.13\textwidth}
		$t = 0.5:$
		\vfill
	\end{minipage}
	\begin{minipage}[]{0.75\textwidth}
		\includegraphics[trim = 50 140 50 130, clip, width=1\textwidth]{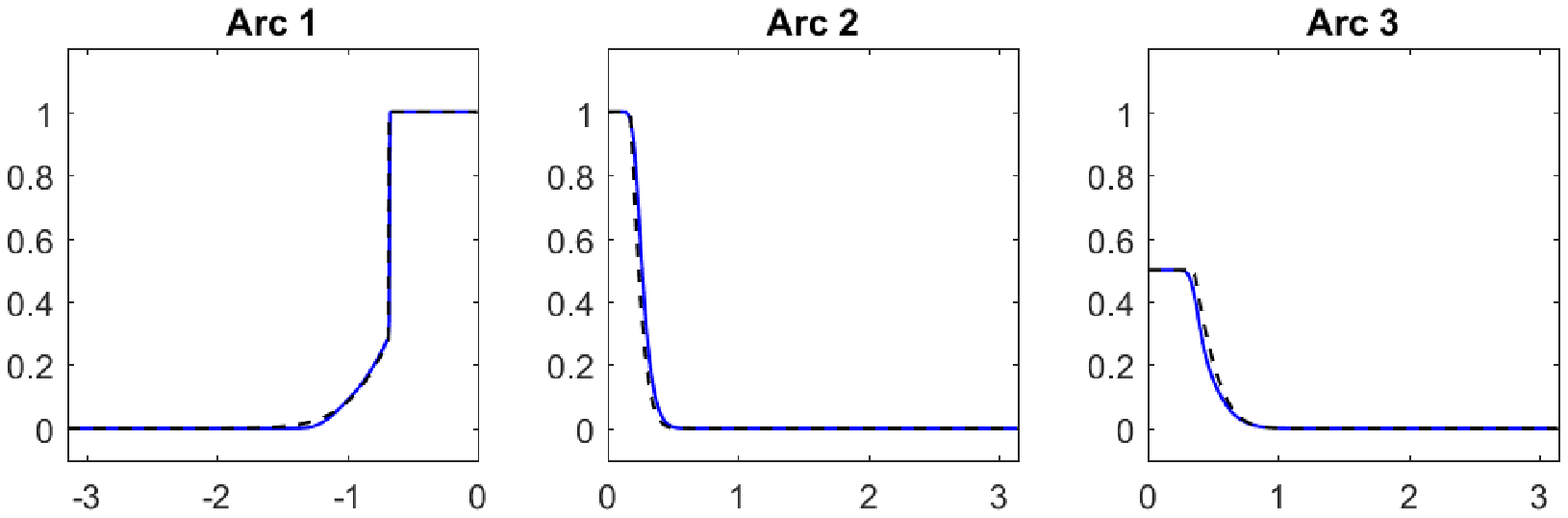}
	\end{minipage}\\
	%	\begin{minipage}[b]{0.13\textwidth}
	%		t = 0.8:
	%		\vfill
	%	\end{minipage}
	%	\begin{minipage}[]{0.77\textwidth}
	%		\includegraphics[trim = 50 140 50 130, clip, width=1\textwidth]{Pictures/Network/112_k08_passive.eps}
	%	\end{minipage}\\
	\begin{minipage}[b]{0.13\textwidth}
		$t = 1.0:$
		\vfill
	\end{minipage}
	\begin{minipage}[]{0.75\textwidth}
		\includegraphics[trim = 50 140 50 130, clip, width=1\textwidth]{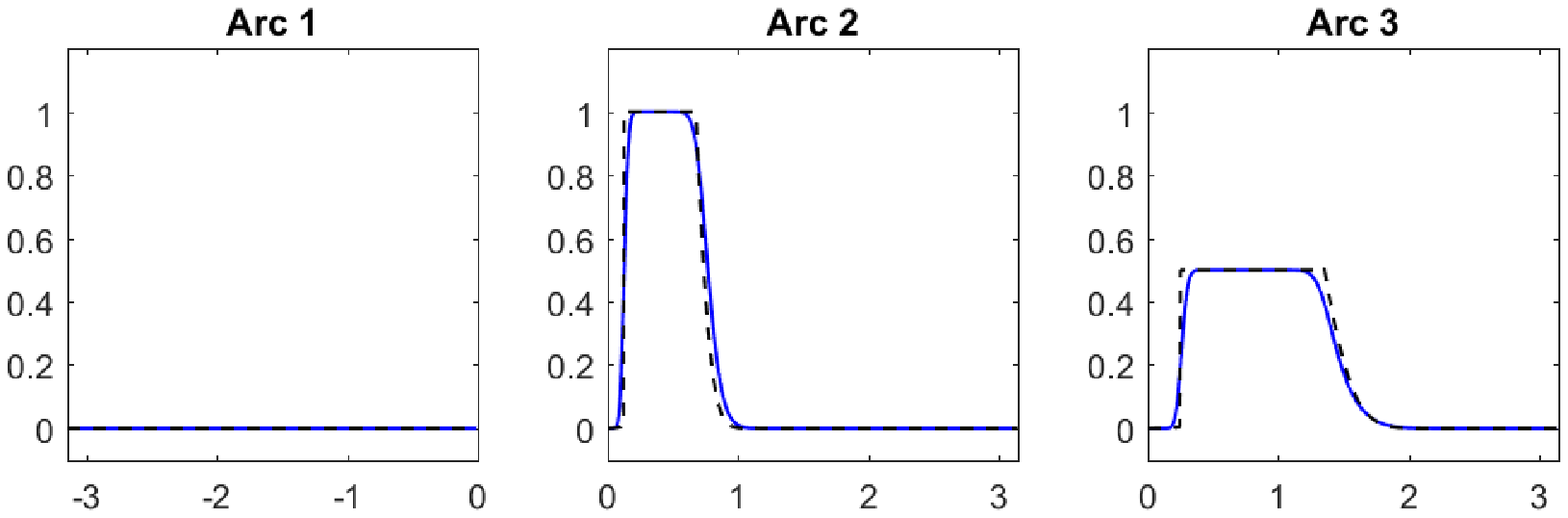}
	\end{minipage}
	\caption{Test 3: ``passive'' junction with distribution parameter $\mu = 0.5$}
	\label{fig:network_test3_passive}
\end{figure}

\subsection{One-to-two junction} 
We pass to the one-to-two junction described in Section \ref{s122}. We consider the domain
$\Omega_1=(-\pi,0)\times \{0\}, ~\{0\}\times\Omega_2=(0,-\pi) , ~\Omega_3=  (0,\pi)\times\{0\}, ~ t\in [0,2]$
with intersection point $x=(0,0)$. We choose the velocities $a\sub{1}=4, a\sub{2}=1, a\sub{3}=2$ and distribution parameter $\mu=0.5$. 
The initial solution $\rhozeronetwork$ is again \eqref{eq:network_test_rhozero}.
As already mentioned, the flux conservation \eqref{eq:network_1to2_junctioncond} is not sufficient to ensure uniqueness of the solution at the intersection
and therefore ``passive'' and ``active'' junctions are considered.
Due to our choice of parameters, there exists a time $t\in[0,2]$ such that the conditions \eqref{regular122} and \eqref{regular122b} are violated and thus 
congestion arises.

\subsubsection*{Test 3: ``passive'' junction}
In the passive junction case, the solution is given by \eqref{sol122_ratio}. 
The comparison of the numerical and analytical solution are displayed in Figure \ref{fig:network_test3_passive}.
In the first column, the density distribution on arc $1$ at different time steps is shown. In the second and the third column, the density distribution on arcs $i = 2, 3$ are presented. 
Congestion starts at about $t_0 = 0.3$, so at time $t = 0.5$, we are already in the congested phase. On arc $i = 3$, a density value of $0.5$ is kept. This is due to the constant ratio of the outgoing fluxes, even if the maximal capacity is not used. As we see here, this leads to shocks in the solution also if the corresponding arc is not congested. 
At time $t = 1.0$, all mass passed the junction and is transported by the outgoing arcs. The final time of congestion is about $t^E = 0.9$ in this scenario.
%We observe an excellent fitting between the approximated and the analytical solution.

\begin{figure}[tb]
	\centering
	\begin{minipage}[b]{0.13\textwidth}
		$t = 0.5:$
		\vfill
	\end{minipage}
	\begin{minipage}[]{0.75\textwidth}
		\includegraphics[trim = 50 140 50 130, clip, width=1\textwidth]{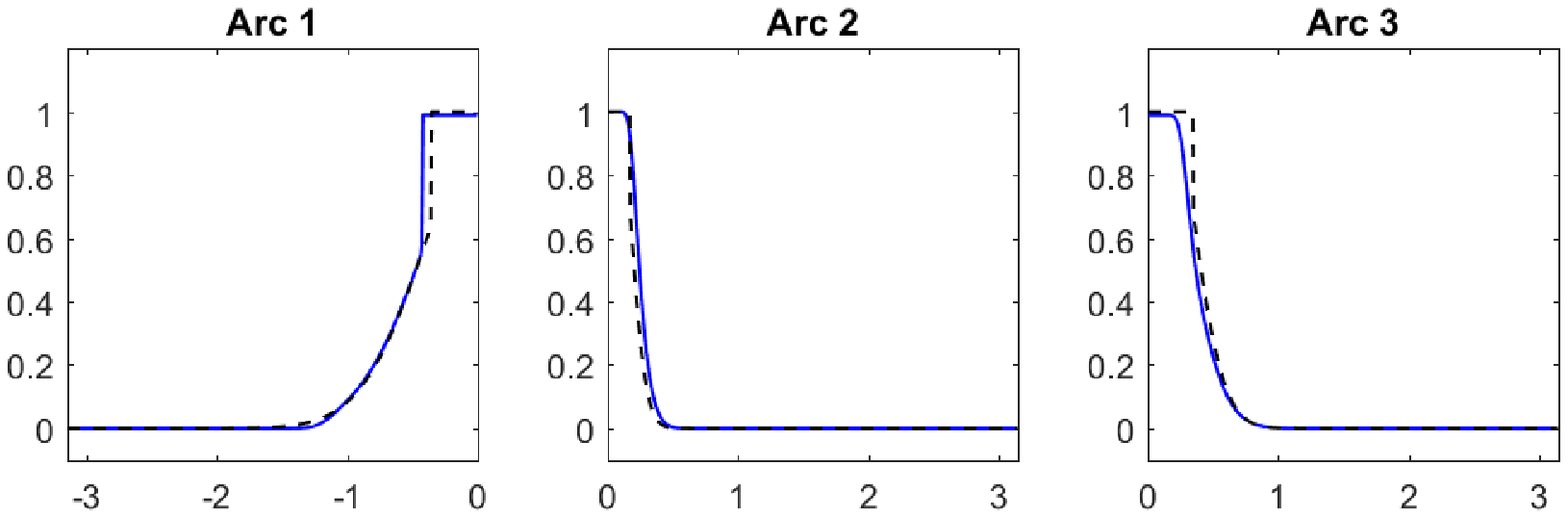}
	\end{minipage}\\
	\begin{minipage}[b]{0.13\textwidth}
		$t = 1.0:$
		\vfill
	\end{minipage}
	\begin{minipage}[]{0.75\textwidth}
		\includegraphics[trim = 50 140 50 130, clip, width=1\textwidth]{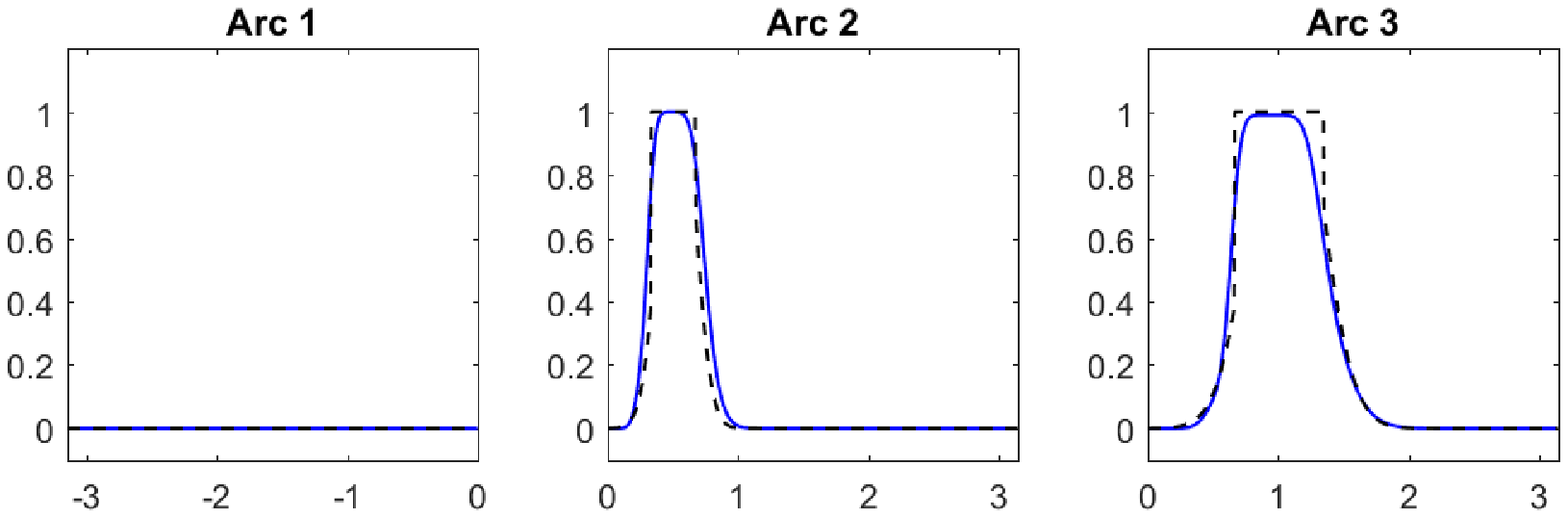}
	\end{minipage}
	\caption{Test 4: ``active'' junction with distribution parameter $\mu = 0.5$}
	\label{fig:network_test4_active}
\end{figure}

\subsubsection*{Test 4: ``active'' junction}
Here, the solution is given by \eqref{sol122_active}. The results in Figure \ref{fig:network_test4_active}
are again for the congested case. 
Note that now congestion starts at about $t_0 = 0.4$.
At time $t = 0.5$, the maximal capacity of both outgoing arcs $i = 2, 3$ is reached and we are in the congested phase. 
Compared to the passive junction case, congestion is reduced on the incoming arc $i = 1$.
At time $t = 1.0$, all mass passed the junction and is transported by the outgoing arcs.
Now, the final time of congestion is about $t^E = 0.7$ which is less than in the previous case. %The overall congestion time is hence reduced. 

\subsection{Two-to-one junction} 
\begin{figure}[tb]
	\centering
	\begin{minipage}[b]{0.1\textwidth}
		$t = 0:$
		\vfill
	\end{minipage}
	\begin{minipage}[]{0.77\textwidth}
		\includegraphics[trim = 50 140 50 130, clip, width=1\textwidth]{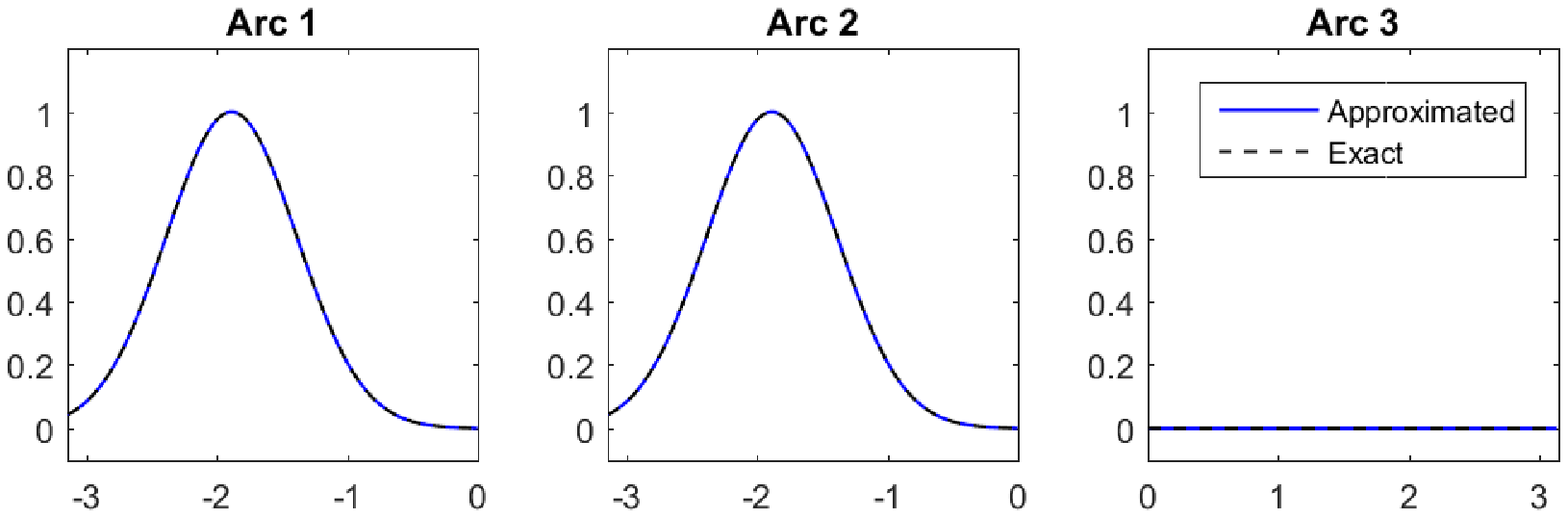}
	\end{minipage}\\
%	\begin{minipage}[b]{0.1\textwidth}
%		t = 1.5:
%		\vfill
%	\end{minipage}
%	\begin{minipage}[]{0.77\textwidth}
%		\includegraphics[trim = 50 140 50 130, clip, width=1\textwidth]{./Figure/221_k15_mu03_fin.eps}
%	\end{minipage}\\
	\begin{minipage}[b]{0.1\textwidth}
		$t = 2:$
		\vfill
	\end{minipage}
	\begin{minipage}[]{0.77\textwidth}
		\includegraphics[trim = 50 140 50 130, clip, width=1\textwidth]{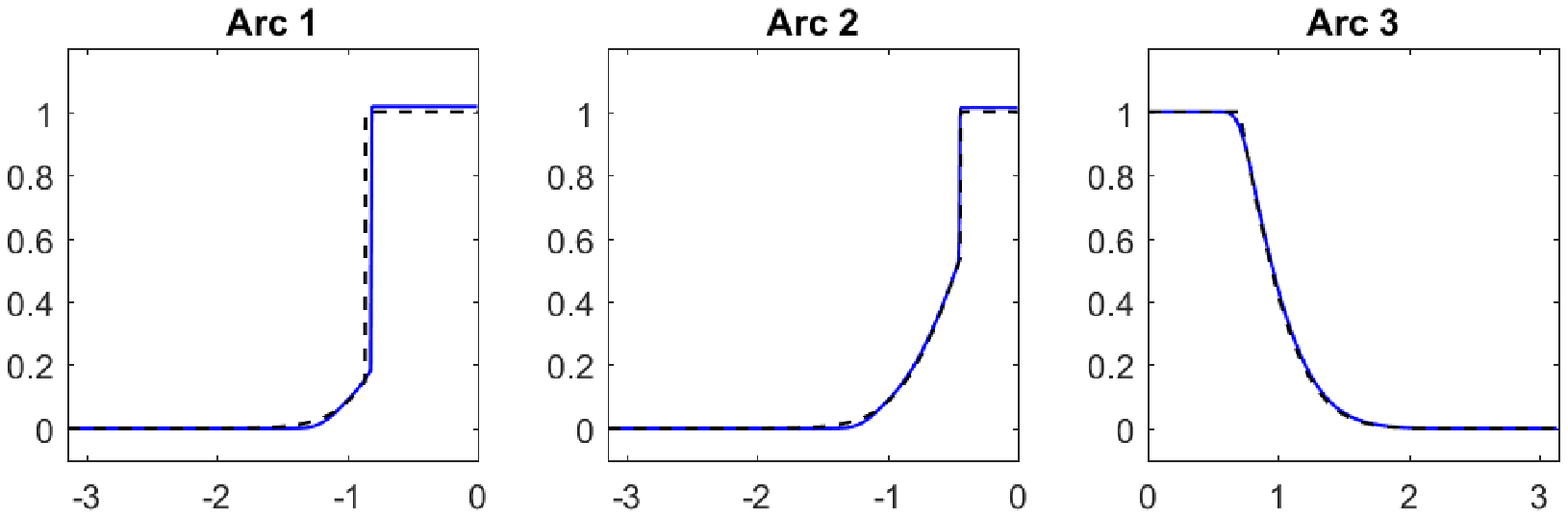}
	\end{minipage}\\
%	\begin{minipage}[b]{0.1\textwidth}
%		t = 3:
%		\vfill
%	\end{minipage}
%	\begin{minipage}[]{0.77\textwidth}
%		\includegraphics[trim = 50 140 50 130, clip, width=1\textwidth]{./Figure/221_k30_mu03_fin.eps}
%	\end{minipage}\\
	\begin{minipage}[b]{0.1\textwidth}
		$t = 4:$
		\vfill
	\end{minipage}
	\begin{minipage}[]{0.77\textwidth}
		\includegraphics[trim = 50 140 50 130, clip, width=1\textwidth]{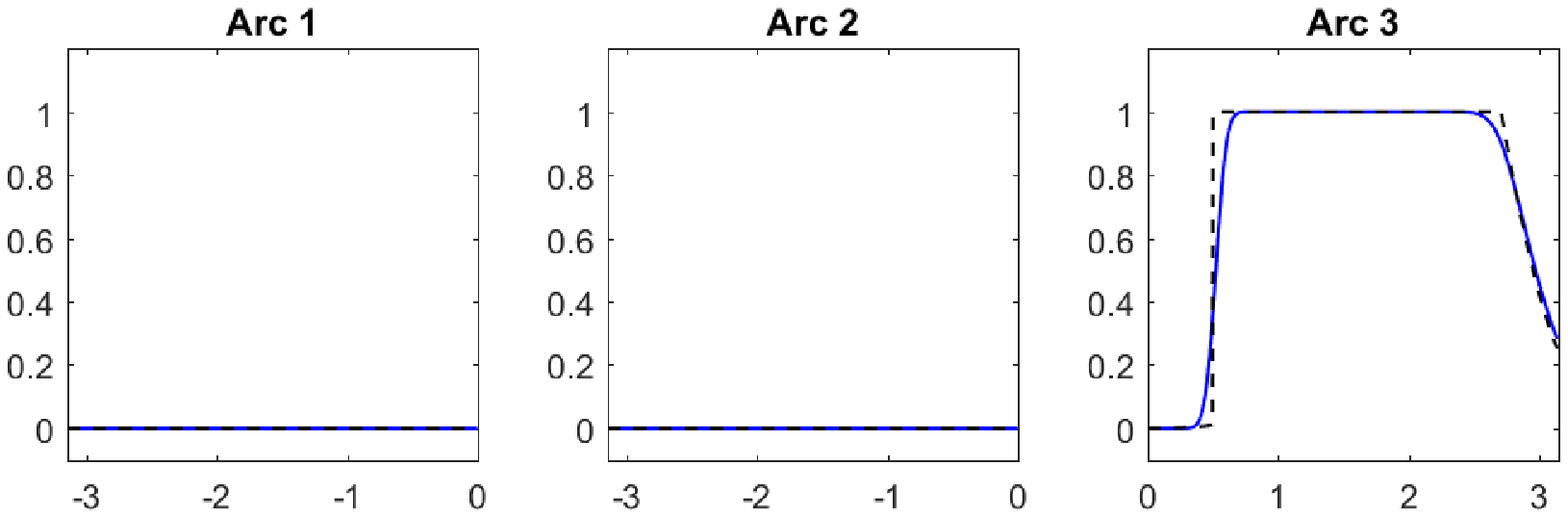}
	\end{minipage}
	\caption{Test 5: merging junction with parameter $q = 0.3$}
	\label{fig:network_221_mu03}
\end{figure}

The last test is the merging junction discussed in Section \ref{s221}. We consider the domain
$\Omega_1=(-\pi,0)\times \{0\}, ~\Omega_2= \{0\}\times(-\pi,0), ~\Omega_3=(0,\pi)\times \{0\},~t\in [0,4].$
The initial data $\rhozeronetwork_i$ on each incoming arc $i = 1, 2$ is \eqref{eq:network_test_rhozero}. On the outgoing arc $i = 3$, 
we set $\rhozeronetwork_3=0.$. 
All velocities are fixed to $a_i=1$ for all arcs. 
%We observe in this scenario a congested behavior despite equal velocities and capacities. 
This setting also allows to recover the results developed in \cite{gottlich2013discontinuous}, where the capacity and the speed are assumed to be equal for all arcs.
Figure \ref{fig:network_221_mu03} shows the result of the evolution of the density at various time steps with merging parameter $q = 0.3$. 
The latter leads to a prioritization of arc $2$ and a non-symmetric transportation on the two incoming arcs.
We observe how the density initially placed on $\Omega_1$ and $\Omega_2$ is transported till $x=(0,0)$ is reached and congestion forms. 
At time $t = 2$, both arcs are congested. Due to the prioritization of arc $2$, congestion is less than on arc $1$.
%At time $t = 3$, on arc 1, we still have a congested area, whereas on arc 2 the congestion vanishes.
At time $t = 4$, all mass is absorbed by the outgoing arc. 
It is worth noticing that, once all density is absorbed by the outgoing arc, the configuration on the outgoing arc is the same, independent of the merging parameter $q$. This is due to the condition \eqref{eq:network_2to1_maximalitycond}, which implies that the outgoing arc always absorbs as much mass as possible.

To conclude, one can say that the numerical method presented in this section is a powerful tool to approximate the most meaningful solution of the material flow problem in the network case given by \eqref{eq:model_network} - \eqref{eq:network_flux_arci}. This is particularly relevant in cases where we can no longer directly compute the analytical solution of the problem. To validate the results, in our numerical simulation study we considered special cases, for which we already derived the analytical solution. For those, the numerical results obtained match very well the analytical solution and catches the behavior of the evolution of the congested area before the junction point. Moreover, an error analysis for the case of a one-to-one junction has shown that the numerical solution converges against the analytical one with respect to the discretized version of the time-averaged $L^2$-norm.

%
%\section{Conclusions}
%
%
%
%
%
%\section*{Acknowledgments}
%We would like to acknowledge ...

\bibliographystyle{siamplain}

\end{document}